\documentclass[12pt,leqno]{article}
\pdfoutput=1 
\usepackage{amssymb}
\usepackage[mathscr]{eucal}
\usepackage{amsmath,amssymb,latexsym,theorem,bbm,graphicx}
\newcommand{\DS}{\displaystyle}

\newcommand{\SC}{\scriptstyle}

\newcommand{\CC}{\mathsf{C}}
\newcommand{\DD}{\mathsf{D}}
\newcommand{\CCb}{\mathbb{C}}

\newcommand{\NN}{\mathbb{N}}

\newcommand{\RR}{\mathbb{R}}

\newcommand{\ZZ}{\mathbb{Z}}

\newcommand{\bA}{{\boldsymbol{A}}}
\newcommand{\btA}{{\boldsymbol{\widetilde A}}}
\newcommand{\bB}{{\boldsymbol{B}}}
\newcommand{\bC}{{\boldsymbol{C}}}
\newcommand{\bD}{{\boldsymbol{D}}}
\newcommand{\be}{{\boldsymbol{e}}}
\newcommand{\bI}{{\boldsymbol{I}}}
\newcommand{\bM}{{\boldsymbol{M}}}
\newcommand{\bP}{{\boldsymbol{P}}}
\newcommand{\bu}{{\boldsymbol{u}}}
\newcommand{\bv}{{\boldsymbol{v}}}

\newcommand{\bX}{{\boldsymbol{X}}}
\newcommand{\bY}{{\boldsymbol{Y}}}

\newcommand{\bPi}{{\boldsymbol{\Pi}}}
\newcommand{\bone}{{\boldsymbol{1}}}

\newcommand{\cA}{{\mathcal A}}

\newcommand{\cD}{{\mathcal D}}

\newcommand{\cF}{{\mathcal F}}

\newcommand{\cL}{{\mathcal L}}
\newcommand{\cM}{{\mathcal M}}

\newcommand{\cU}{{\mathcal U}}
\newcommand{\cX}{{\mathcal X}}
\newcommand{\cY}{{\mathcal Y}}
\newcommand{\cW}{{\mathcal W}}

\newcommand{\dd}{\mathrm{d}}

\newcommand{\slu}{{\SC\mathrm{lu}}}

\newcommand{\INARp}{\textup{INAR($p$)}}

\newcommand{\EE}{\operatorname{E}}
\newcommand{\PP}{\operatorname{P}}
\newcommand{\var}{\operatorname{Var}}

\newcommand{\SE}{\operatorname{SE}}

\newcommand{\hM}{\widehat{M}}

\newcommand{\hmue}{\widehat{\mu_\vare}}

\newcommand{\halpha}{\widehat{\alpha}}

\newcommand{\tS}{\widetilde{S}}

\newcommand{\varr}{\varrho}
\newcommand{\vare}{\varepsilon}
\newcommand{\tvare}{\widetilde{\vare}}

\renewcommand{\mid}{\,|\,}
\newcommand{\bmid}{\,\big|\,}

\renewcommand{\leq}{\leqslant}
\renewcommand{\geq}{\geqslant}

\newcommand{\stoch}{\stackrel{\PP}{\longrightarrow}}
\newcommand{\distr}{\stackrel{\cL}{\longrightarrow}}

\newcommand{\lu}{\stackrel{\slu}{\longrightarrow}}

\newcommand{\bbone}{\mathbbm{1}}
\newcommand{\ns}{{\lfloor ns\rfloor}}
\newcommand{\nt}{{\lfloor nt\rfloor}}
\newcommand{\nT}{{\lfloor nT\rfloor}}
\newcommand{\proofend}{\hfill\mbox{$\Box$}}

\numberwithin{equation}{section}

\theoremstyle{change} \theorembodyfont{\em}
\newtheorem{Lem}{Lemma.}[section]
\newtheorem{Thm}{Theorem.}[section]
\newtheorem{Pro}{Proposition.}[section]
\newtheorem{Cor}{Corollary.}[section]
\newtheorem{Def}{Definition.}[section]

\theorembodyfont{\rm}
\newtheorem{Rem}{Remark.}[section]

\begin{document}

\begin{center}
 {\bfseries\Large Asymptotic behavior of unstable \INARp\ processes} \\[5mm]

 {\sc\large M\'aty\'as $\text{Barczy}^{*,\diamond}$,
            \ M\'arton $\text{Isp\'any}^*$, \ Gyula $\text{Pap}^{\sharp}$}
\end{center}

\vskip0.2cm

\noindent * Faculty of Informatics, University of Debrecen,
            Pf.~12, H--4010 Debrecen, Hungary.

\noindent $\sharp$ Bolyai Institute, University of Szeged, Aradi v\'ertan\'uk tere 1,
                   H--6720 Szeged, Hungary.


\noindent e--mails: barczy.matyas@inf.unideb.hu (M. Barczy),
                    ispany.marton@inf.unideb.hu (M. Isp\'any),
                    papgy@math.u-szeged.hu (G. Pap).

\noindent $\diamond$ Corresponding author.



\renewcommand{\thefootnote}{}
\footnote{\textit{2000 Mathematics Subject Classifications\/}:
          60J80,  60F99.}
\footnote{\textit{Key words and phrases\/}:
 unstable \INARp\ processes, squared Bessel processes, Boston armed robberies data set}

\vspace*{-10mm}

\begin{abstract}
In this paper the asymptotic behavior of an unstable integer-valued
 autoregressive model of order $p$ (\INARp) is described.
Under a natural assumption it is proved that the sequence of appropriately
 scaled random step functions formed from an unstable \INARp\ process converges
 weakly towards a squared Bessel process.
We note that this limit behavior is quite different from that of familiar unstable autoregressive
 processes of order $p$.
An application for Boston armed robberies data set is presented.
\end{abstract}

\section{Introduction}

Recently, there has been remarkable interest in integer-valued
 time series models and a number of results are now available
 in specialized monographs (e.g., MacDonald and Zucchini \cite{MacZuc},
 Cameron and Trivedi \cite{CamTri}, and Steutel and van Harn \cite{SteHar2})
 and review papers (e.g., McKenzie \cite{McK2},
 Jung and Tremayne \cite{JunTre}, and Wei{\ss} \cite{Wei1}).
Reasons to introduce discrete data models come from the need to account for
 the discrete nature of certain data sets, often counts of events, objects or
 individuals.
Examples of applications can be found in the analysis of time series of count
 data on the area of financial mathematics by analyzing stock transactions
(Quoreshi \cite{Quo}), insurance by modeling claim counts (Gouri\' eroux
and Jasiak \cite{GoJa}), medicine by investigating disease incidence (Cardinal
et al. \cite{CaRoLa}), neurobiology by change-point analysis of neuron
spike train data (B\' elisle et al. \cite{BJMWB}), optimal alarm systems
(Monteiro et al. \cite{MonPerSco}), psychometrics by treating longitudinal
count data (B\" ockenholt \cite{Bo1}, \cite{Bo2}), environmetrics by analyzing rainfall
measurements (Thyregod et al. \cite{TCMA}), experimental biology (Zhou and
Basawa \cite{ZhoBas}), and queueing systems (Ahn et al. \cite{AhnGyeJon}
and Pickands III and Stine \cite{PiSt}).

Among the most successful integer-valued time series models proposed in the
literature we mention the INteger-valued AutoRegressive model of order $p$
\ (\INARp). This model was first introduced by McKenzie \cite{McK}
 and Al-Osh and Alzaid \cite{AloAlz1} for the case $p=1$.
The INAR(1) model has been investigated by several authors.
Franke and Seligmann \cite{FrSe} analyzed maximum likelihood estimation
of parameters under Poisson innovation. Du and Li \cite{DuLi} and
Freeland and McCabe \cite{FrMC} derived the limit-distribution of the
conditional least squares estimator of the autoregressive parameter. Silva and Oliveira \cite{SilOli1}
proposed a frequency domain based estimator, Br\"{a}nn\"{a}s and Hellstr\"{o}m
\cite{BraHel} investigated generalized method of moment estimation, Silva and Silva \cite{SilSil} considered
a Yule-Walker estimator. Jung et al. \cite{JRT} analyzed the finite sample
behavior of several estimators  by a Monte Carlo study. Isp\' any et al.
\cite{IspPapZui1}, \cite{IspPapZui2} derived asymptotic inference for nearly
unstable INAR(1) models which has been refined by Drost et al. \cite{DroAkkWer2}
later. A Poisson limit theorem has been proved for an inhomogeneous nearly
critical INAR(1) model by Gy\" orfi et al. \cite{GyorIspPapVar}.

The more general \INARp\ processes were first introduced by Al-Osh and Alzaid
 \cite{AloAlz2}.
In their setup the autocorrelation structure of the process corresponds to that of
 an ARMA($p,p-1$) process, see also Section \ref{INARp_section}.
Another definition of an INAR(p) process was proposed independently
 by Du and Li \cite{DuLi} and by Gauthier and Latour \cite{GauLat} and Latour \cite{Lat2},
 and is different from that of Alzaid and Al-Osh \cite{AloAlz2}.
In Du and Li's setup the autocorrelation structure
 of an \INARp\ process is the same as that of an AR($p$) process. The setup of Du and Li
 \cite{DuLi} has been followed by most of the authors, and our approach will also be
 the same, see Section \ref{INARp_section}.
The \INARp\ model has been investigated by several authors from different points of views.
 Drost et al. \cite{DroAkkWer1} provided asymptotically efficient
 estimator for the parameters.
Silva and Oliveira \cite{SilOli2} described the higher order moments and cumulants
 of \INARp\ processes, and Silva and Silva \cite{SilSil} derived asymptotic distribution of
 the Yule-Walker estimator.
Drost et al. \cite{DroAkkWer3} considered semiparametric INAR($p$) models and proposed efficient
 estimation for the autoregression parameters and innovation distributions.
Recently, the so called $p$-order Rounded INteger-valued AutoRegressive (RINAR($p$)) time series model
 was introduced and studied by Kachour and Yao \cite{KacYao} and Kachour \cite{Kac}.
The broad scope of the empirical literature in which INAR models are applied indicates its relevance.
Examples of such applications include Franke and Seligmann \cite{FrSe} (epileptic seizure counts),
 B\"ockenholt \cite{Bo2} (longitudional count data), Thyregod et al. \cite{TCMA} (rainfall measurements),
 Br\"{a}nn\"{a}s and Hellstr\"{o}m  \cite{BraHel} and Rudholm \cite{Rud} (economics),
 Br\"{a}nn\"{a}s and Shahiduzzaman \cite{BraSha} (finance), Gourieroux and Jasiak \cite{GoJa} (insurance),
 Pavlopoulos and Karlis \cite{PavKar} (environmental studies)
 and McCabe et al. \cite{McCMarHar} (finance and mortality).

An interesting problem, which has not yet been addressed for INAR($p$)
 models, is to investigate the asymptotic behavior of unstable INAR($p$)
 processes, i.e., when the characteristic polynomial has a unit root.
In this paper we give a complete description of this limit behavior.
In particular, it will turn out that an  INAR($p$) model is unstable if
 and only if the sum of its autoregressive parameters equals 1, and in
 this case the only unit root of the characteristic polynomial is 1 with
 multiplicity one.
For the sake of convenience, we suppose that the process starts from zero.
Without loss of generality, we may suppose that the $p$th autoregressive parameter is strictly positive and
 that the greatest common divisor of the strictly positive autoregressive parameters is 1,
 see Remark \ref{REMARK3}.
Under the assumption that the second moment of the innovation distribution is finite, we prove that
 the sequence of appropriately scaled random step functions formed from
 an unstable \INARp\ process converges weakly towards a squared Bessel process.
This limit process is a continuous branching process  also known as square-root process or
 Cox-Ingersoll-Ross process.
We remark that a similar theorem holds for critical, i.e., when the offspring mean is equal
 to 1, branching processes with immigration, see Wei and Winnicki \cite[Theorem 2.1]{WW1}.
We should also note that the asymptotic behavior of unstable \INARp\ models
 is completely different from that of familiar (real-valued) unstable AR($p$) models in at least
 two senses.
On the one hand, the characteristic polynomial of a primitive (see Definition \ref{DEF_primitive})
 unstable \INARp\ model has only one unit root, namely 1, with multiplicity one,
 whereas for an unstable AR($p$) model it may have real or complex unit roots with various different multiplicities.
On the other hand, in the case of a primitive unstable \INARp\ model there is a limit process which is a squared
 Bessel process, while in the case of an unstable AR($p$) model in general there is no limit process,
 only for appropriately transformed and scaled random step functions, see
 Chan and Wei \cite{CW}, Jeganathan \cite{Jeg} and van der Meer et al. \cite[Theorem 3]{MPZ}.

We remark that our result can be considered as the first step towards the comprehensive theory of
 nonstationary integer-valued time series and investigation of the unit root problem of
 econometrics in the integer-valued setup.
Nonstationary time series have been playing an important role in both econometric
theory and applications over the last 20 years, and a substantial literature
has been developed in this field. A detailed set of references is given in Phillips
and Xiao \cite{PX}.
We note that Ling and Li \cite{LL1}, \cite{LL2} considered an unstable ARMA model
 with GARCH errors and an unstable fractionally integrated ARMA model.
Concerning relevance and practical applications of unstable INAR models we note that
 empirical studies show importance of these kind of models.
Br\"{a}nn\"{a}s and Hellstr\"{o}m  \cite{BraHel} reported an INAR$(1)$ model with a coefficient \ $0.98$
 for the number of private schools, Rudholm \cite{Rud} considered INAR$(1)$ \ models
 with coefficients \ $0.98$ \ and \ $0.99$, respectively for the number of Swedish
 generic-pharmaceutical market.
Hellstr\"{o}m \cite{Hel} focused on the testing of unit root in INAR(1) models and provided
 small sample distributions for the Dickey-Fuller test statistic under the null hypothesis
 of unit root in an INAR(1) model with Poisson distributed innovations.
In this paper, we report that a subset INAR$(12)$ model is an adequate model for
 Boston armed robberies data set published in Deutsch and Alt \cite{DA}.
Our proposed model can be considered unstable since
 the sum of the estimated (autoregressive) coefficients is 1.0317.
To our knowledge a unit root test for general INAR$(p)$ models is not known, and from
 this point of view studying unstable INAR$(p)$ models is an important preliminary task.

The rest of the paper is organized as follows.
Section \ref{INARp_section} provides a background description of
 basic theoretical results related with \INARp\ models.
In Section \ref{conv_section} we describe the asymptotic behavior of unstable \INARp\
 processes.
Under the assumption that the second moment of the innovation distribution is finite,
 we prove that the sequence of appropriately scaled random step functions formed
 from an unstable \INARp\ process converges weakly towards a squared Bessel process,
 see Theorem \ref{main}.
Section \ref{Boston_armrob} presents a real-life application of unstable
 INAR$(p)$ models by investigating the Boston armed robberies time series.
 Section \ref{proof_section} contains a proof of our main Theorem \ref{main}.
For the proof, we collect some properties of the first and second moments of
 (not necessarily unstable) \INARp\ processes, we recall a useful functional martingale
 limit theorem and an appropriate version of the continuous mapping theorem,
 see Lemma \ref{Moments}, Corollary \ref{EEX}, Theorem \ref{Conv2DiffCor}
 and Lemma \ref{Conv2Funct} in Appendix, respectively.

\section{The \INARp\ model}\label{INARp_section}

Let \ $\ZZ_+$, \ $\NN$, \ $\RR$, \ $\RR_+$ \ and \ $\CCb$ \ denote the set of non-negative integers,
 positive integers, real numbers, non-negative real numbers and complex numbers, respectively.
For all \ $n\in\NN$, \ let us denote by \ $\bI_{\mathbf n}$ \ the $n\times n$ identity matrix.
Every random variable will be defined on a fixed probability space
 \ $(\Omega,\cA,\PP)$.

One way to obtain models for integer-valued data is to replace multiplication
 in the conventional ARMA models in such a way to ensure the integer
 discreteness of the process and to adopt the terms of self-decomposability and
 stability for integer-valued time series.

\begin{Def}
Let \ $(\vare_k)_{k\in\NN}$ \ be an independent and identically distributed
 (i.i.d.) sequence of non-negative integer-valued random variables, and let
 \ $\alpha_1,\ldots,\alpha_p \in [0,1]$.
\ An \INARp\ time series model with coefficients \ $\alpha_1,\ldots,\alpha_p$
 \ and innovations \ $(\vare_k)_{k\in\NN}$ \ is a stochastic process
 \ $(X_n)_{n\geq -p+1}$ \ given by
 \begin{align}\label{INARp}
   X_k = \sum_{j=1}^{X_{k-1}} \xi_{k,1,j}
         + \cdots + \sum_{j=1}^{X_{k-p}} \xi_{k,p,j} + \vare_k , \qquad k \in \NN ,
 \end{align}
 where for all \ $k\in\NN$ \ and \ $i\in\{1,\dots,p\}$, \ $(\xi_{k,i,j})_{j\in\NN}$
 \ is a sequence of i.i.d.\ Bernoulli random variables with mean \ $\alpha_i$
 \ such that these sequences are mutually independent and independent of the
 sequence \ $(\vare_k)_{k\in\NN}$, \ and
 \ $X_0$, $X_{-1}$, \ldots, $X_{-p+1}$ \ are non-negative integer-valued random
 variables independent of the sequences \ $(\xi_{k,i,j})_{j\in\NN}$, \ $k\in\NN$,
 \ $i\in\{1,\dots,p\}$, \ and \ $(\vare_k)_{k\in\NN}$.
\end{Def}

The \INARp\ model \eqref{INARp} can be written in another way using
 the binomial thinning operator \ $\alpha\,\circ$
 \ (due to Steutel and van Harn~\cite{SteHar}) which we recall now.
Let \ $X$ \ be a non-negative integer-valued random variable.
Let \ $(\xi_j)_{j\in\NN}$ \ be a sequence of i.i.d.\ Bernoulli random variables
 with mean \ $\alpha\in[0,1]$.
\ We assume that the sequence \ $(\xi_j)_{j\in\NN}$ \ is independent of \ $X$.
\ The non-negative integer-valued random variable \ $\alpha\,\circ X$
 \ is defined by
 \[
   \alpha\circ X
     :=\begin{cases}
        \sum\limits_{j=1}^X\xi_j, & \quad \text{if \ $X>0$},\\[2mm]
         0, & \quad \text{if \ $X=0$}.
       \end{cases}
 \]
The sequence \ $(\xi_j)_{j\in\NN}$ \ is called a counting sequence.
The \INARp\ model \eqref{INARp} takes the form
 \[
    X_k = \alpha_1 \circ X_{k-1} + \cdots + \alpha_p \circ X_{k-p} + \vare_k ,
    \qquad k \in \NN .
 \]
Note that the above form of the \INARp \ model is quite analogous with a usual AR($p$)
 process (another slight link between them is the similarity of some
 conditional expectations, see \eqref{AR_INAR}).
As we noted in the introduction, this definition of the INAR(p) process
 was proposed independently by Du and Li \cite{DuLi} and by Gauthier and Latour \cite{GauLat}
 and Latour \cite{Lat2},
 and is different from that of Alzaid and Al-Osh \cite{AloAlz2}, which assumes that the
 conditional distribution of the vector
 \ $(\alpha_1\circ X_t, \alpha_2\circ X_t, \ldots , \alpha_p\circ X_t)$ \ given \ $X_t = x_t$
 \ is multinomial with parameters \ $(\alpha_1, \alpha_2, \ldots , \alpha_p, x_t)$ \ and
 is independent of the past history of the process.
The two different formulations imply different second-order structure for the processes:
 under the first approach, the \INARp\ has the same second-order structure as an AR($p$)
 process, whereas under the second one, it has the same one as an ARMA$(p, p-1)$
 process.

An alternative representation of the \INARp \ process as a p-dimensional INAR(1)
 process was obtained by Franke and Subba Rao \cite{FraSub} and see also Latour \cite[formula (2.3)]{Lat1}.
Accordingly, the \INARp \ process defined in \eqref{INARp} can be written as
 \[
    \bX_k = \bA \boldsymbol{\circ} \bX_{k-1} + \boldsymbol{\vare}_k, \qquad k\in\NN,
 \]
 where the \ $p$-dimensional random vectors \ $\bX_k$, $\boldsymbol{\vare}_k$ \ and the \ $(p\times p)$-matrix
 \ $\bA$ \ are defined by
 \begin{align}\label{HELP_MATRIX_IRASMOD}
   \bX_k := \begin{bmatrix}
             X_k \\ X_{k-1} \\ X_{k-2} \\ \vdots \\ X_{k-p+2} \\ X_{k-p+1}
            \end{bmatrix} , \qquad
   \bA := \begin{bmatrix}
            \alpha_1 & \alpha_2 & \alpha_3 & \cdots & \alpha_{p-1} & \alpha_p \\
               1    &     0   &     0    & \cdots &      0      &     0    \\
               0    &     1   &     0    & \cdots &      0      &     0    \\
           \vdots   & \vdots  & \vdots  &  \ddots & \vdots      & \vdots   \\
               0    &     0   &     0    & \cdots &      0      &     0    \\
               0    &     0   &     0    & \cdots &      1      &     0
          \end{bmatrix} , \qquad
   \boldsymbol{\vare}_k
   := \begin{bmatrix} \vare_k \\ 0 \\ 0 \\ \vdots \\ 0 \\ 0 \end{bmatrix},
 \end{align}
 and for a $p$-dimensional random vector \ $\bY=(Y_1,\ldots,Y_p)$ \ and
 for a \ $p\times p$ \ matrix \ $\bB=(b_{ij})_{i,j=1}^p$
 \ with entries satisfying \ $0\leq b_{ij}\leq 1$, $i,j=1,\ldots,p$, \ the matricial binomial thinning operation
 \ $\bB\boldsymbol{\circ} \bY$ \ is defined as a \ $p$-dimensional random vector whose \ $i$-th component,
 \ $i=1,\ldots,p$, \ is given by
 \[
    \sum_{j=1}^p b_{ij}\circ Y_j,
 \]
 where the counting sequences of all \ $b_{ij}\circ Y_j$, $i,j=1,\ldots,p$, \ are assumed to be independent
 of each other.

In what follows for the sake of simplicity we consider a zero start \INARp\ process,
 that is we suppose \ $X_0 = X_{-1} = \ldots = X_{-p+1} = 0$.
The general case of nonzero initial values may be handled in a similar way, but we renounce
 to consider it.
For nonzero initial values the first and second order moments of the sequence
 \ $(X_k)_{k\in\ZZ_+}$ \ have a more complicated form than in Lemma \ref{Moments}.
Further, for proving a corresponding version of our main result (see Theorem \ref{main})
 one needs to apply a more general version of Theorem \ref{Conv2DiffCor} which is also valid
 for random step functions not necessarily starting from $0$.

In the sequel, we always assume that \ $\EE(\vare_1^2) < \infty$.
\ Let us denote the mean and variance of \ $\vare_1$ \ by \ $\mu_\vare$ \ and
 \ $\sigma_\vare^2$, \ respectively.

For all \ $k\in\ZZ_+$, \ let us denote by \ $\cF_k$ \ the \ $\sigma$-algebra generated
 by the random variables \ $X_0,X_1,\ldots,X_k$. (Note that \ $\cF_0=\{\emptyset,\Omega\}$, \
 since \ $X_0=0$.)
\ By \eqref{INARp},
 \begin{equation}\label{AR_INAR}
   \EE(X_k\mid\cF_{k-1})
   = \alpha_1 X_{k-1} + \cdots + \alpha_p X_{k-p} + \mu_\vare , \qquad k \in \NN.
 \end{equation}
Consequently,
 \[
   \EE(X_k)
   = \alpha_1 \EE(X_{k-1}) + \cdots + \alpha_p \EE(X_{k-p}) + \mu_\vare ,
   \qquad k \in \NN.
 \]
This can also be written in the form
 \ $\EE(\bX_k) = \bA \EE(\bX_{k-1}) + \mu_\vare \be_1$, \ $k \in \NN$, \ where
 \ $
   \be_1
   := \begin{bmatrix} 1 , 0 , 0 , \ldots , 0 , 0 \end{bmatrix}^\top
     \in\RR^{p\times 1}.
   $
\ Consequently, we have
 \[
   \EE(\bX_k) = \mu_\vare \sum_{j=0}^{k-1} \bA^j \be_1 , \quad k\in\NN,
 \]
 which implies
 \begin{equation}\label{expect_rec}
   \EE(X_k) = \EE(\be_1^\top \bX_k)
            = \mu_\vare \sum_{j=0}^{k-1} \be_1^\top \bA^j \be_1 , \quad k\in\NN.
 \end{equation}
Hence the matrix \ $\bA$ \ plays a crucial role in the description of asymptotic behavior
 of the sequence \ $( X_k )_{ k\geq -p+1}$.
\ Let \ $\varr(\bA)$ \ denote the spectral radius of \ $\bA$, \ i.e.,
 the maximum of the modulus of the eigenvalues of \ $\bA$.

In what follows we collect some known facts about the matrix \ $\bA$.
\ First we recall the notions of irreducibility and primitivity of a matrix.
A matrix \ $\bM \in \RR^{p \times p}$ \ is called reducible if \ $p = 1$ \ and
 \ $\bM = 0$, \ or if \ $p \geq 2$ \ and there exist a permutation matrix
 \ $\bP \in \RR^{ p \times p}$ \ and an integer \ $r$ \ with
 \ $1 \leq r \leq p-1$ \ such that
 \[
  \bP^\top \bM \bP = \begin{bmatrix} \bB & \bC \\ {\boldsymbol{0}} & \bD \end{bmatrix},
 \]
 where \ $\bB \in \RR^{ r \times r}$, \ $\bD \in \RR^{ (p-r) \times (p-r) }$,
 \ $\bC \in \RR^{ r \times (p-r) }$, \ and \ ${\boldsymbol{0}} \in \RR^{ (p-r) \times r}$
 \ is a null matrix.
A matrix \ $\bM \in \RR^{ p \times p}$ \ is called irreducible if it is not
 reducible, see, e.g., Horn and Johnson \cite[Definitions 6.2.21 and 6.2.22]{HJ}.
A matrix \ $\bM \in \RR_+^{ p \times p}$ \ is called primitive if it is irreducible
 and has only one eigenvalue of maximum modulus, see, e.g., Horn and Johnson \cite[Definition 8.5.0]{HJ}.
By Horn and Johnson \cite[Theorem 8.5.2]{HJ}, a matrix \ $\bM \in \RR_+^{ p \times p}$ \ is primitive
 if and only if there exists a positive integer \ $k$ \ such that all the entries of
 the matrix \ $\bM^k$ \ are positive.

Let us denote by \ $\varphi$ \ the characteristic polynomial of the matrix \ $\bA$, \ i.e.,
 \[
   \varphi(\lambda)
   := \det(\lambda \bI_{\mathbf p} - \bA)
   = \lambda^p - \alpha_1 \lambda^{p-1} - \cdots - \alpha_{p-1} \lambda - \alpha_p,
   \qquad \lambda\in \CCb.
 \]

\begin{Pro}\label{Proposition2}
For \ $\alpha_1,\ldots,\alpha_p\in[0,1]$, \ $\alpha_p>0$, \ let us consider the matrix \ $\bA$ \ defined
 in \eqref{HELP_MATRIX_IRASMOD}.
Then the following assertions hold:
 \renewcommand{\labelenumi}{{\rm(\roman{enumi})}}
 \begin{enumerate}
   \item The characteristic polynomial \ $\varphi$ \ has just one positive root,
         \ $\varr(\bA) > 0$, \ the nonnegative matrix \ $\bA$ \ is irreducible,
         \ $\varr(\bA)$ \ is an eigenvalue of \ $\bA$ \ and
         \begin{align} \label{char_poli0}
            & \sum_{k=1}^p \alpha_k \varrho(\bA)^{-k} = 1, \\ \label{char_poli1}
            & \sum_{k=1}^p k \alpha_k \varrho(\bA)^{-k} = \varrho(\bA)^{-p+1}\varphi^\prime(\varrho(\bA)).
         \end{align}
        Further,
        \begin{align}\label{char_poli2}
          & \varrho(\bA)\,
             \begin{cases}
              < & \\
              = & \\
              >
             \end{cases}
              1
              \qquad \Longleftrightarrow  \qquad
            \sum_{k=1}^p\alpha_k\,
               \begin{cases}
              < & \\
              = & \\
              >
             \end{cases}
              1.
        \end{align}
   \item If the greatest common divisor \ $d$ \ of the set
          \ $\big\{i\in\{1,\ldots,p\} : \alpha_i>0 \big\}$ is equal to one, then \ $\bA$ \ is primitive,
          \ $\varr(\bA)$ \ is an eigenvalue of \ $\bA$, \ the algebraic and geometric multiplicity of
          \ $\varr(\bA)$ \ equal 1 and the absolute value of the other eigenvalues of \ $\bA$ \ are
          less than \ $\varr(\bA)$.
          \ Corresponding to the eigenvalue \ $\varr(\bA)$ \ there exists a unique vector
          \ $\bu_\bA \in \RR^p$ \ with positive coordinates such that
          \ $\bA \bu_\bA = \varrho(\bA) \bu_\bA$ \ and the sum of the coordinates of \ $\bu_A$ \ is 1,
          namely, \ $\bu_\bA$ \ takes the form
          \[
             \bu_\bA = \begin{bmatrix}
                         \bu_{\bA,1} , \ldots , \bu_{\bA,p} \\
                       \end{bmatrix}^\top
               \qquad \text{with} \qquad
                   \bu_{\bA,i} :=  \frac{\varrho(\bA)^{-i+1}}
                           {\sum_{k=1}^p \varrho(\bA)^{-k+1}},\quad i=1,\ldots,p.
          \]
          Further,
          \begin{equation}\label{conv}
             \varrho(\bA)^{-n} \bA^n \to \bPi_\bA := \bu_\bA \bv_\bA^\top ,
              \qquad \text{as \ $n\to\infty$,}
          \end{equation}
          where \ $\bv_\bA \in \RR^p$ \ is a unique vector with positive coordinates
          such that \ $\bA^\top \bv_\bA = \varrho(\bA) \bv_\bA$ \ and
          \ $\bu_\bA^\top \bv_\bA = 1$, \ namely \ $\bv_\bA$ \ takes the form
          \ $
           \bv_\bA  = \begin{bmatrix}
                         \bv_{\bA,1} , \ldots , \bv_{\bA,p} \\
                       \end{bmatrix}^\top
          $ \ with
          \begin{align*}
                   \bv_{\bA,i}
                       := \frac{\sum_{k=1}^p \varrho(\bA)^{-k+1}}
                                   {\sum_{k=1}^p k\alpha_k\varrho(\bA)^{-k}}
                                 \sum_{\ell=i}^p \alpha_\ell \varrho(\bA)^{i-1-\ell}
                       = \frac{\sum_{k=1}^p \varrho(\bA)^{-k+1}}
                                 {\varrho(\bA)^{-p+1}\varphi^\prime(\varrho(\bA))}
                                 \sum_{\ell=i}^p \alpha_\ell \varrho(\bA)^{i-1-\ell},
          \end{align*}
          for \ $i=1,\ldots,p$.
          \ Moreover, there exist positive numbers \ $c_\bA$ \ and \ $r_\bA$ \ with
          \ $r_\bA < 1$ \ such that for all \ $n \in \NN$
          \begin{equation}\label{rate}
                 \| \varrho(\bA)^{-n} \bA^n - \bPi_\bA \| \leq c_\bA r_\bA^n ,
          \end{equation}
          where \ $\| \bB \|$ \ denotes the operator norm of a matrix
          \ $\bB \in \RR^{p \times p}$ \ defined by
          \ $\| \bB \| := \sup_{ \| x \| =1 } \| \bB x \|$.
 \end{enumerate}
\end{Pro}

\noindent
\textbf{Proof.}
{(i):} First we check that \ $\varphi$ \ has just one positive root, which readily yields that
 \ $\varrho(\bA)>0$.
\ The function
 \ $\lambda
    \mapsto 1 - \lambda^{-p} \varphi(\lambda)
    = \alpha_1 \lambda^{-1} + \cdots + \alpha_{p-1} \lambda^{-p+1}
                            + \alpha_p \lambda^{-p}$
 \ is strictly decreasing and continuous on \ $(0,\infty)$ \ with
 \ $\lim\limits_{\lambda \downarrow 0} (1 - \lambda^{-p} \varphi(\lambda)) = \infty$
 \ and
 \ $\lim\limits_{\lambda \uparrow \infty} (1 - \lambda^{-p} \varphi(\lambda)) = 0$,
 \ thus it takes the value 1 at exactly one positive point, which is the only
 positive root of \ $\varphi$.

\noindent Now we turn to check that \ $\bA$ \ is irreducible.
By Brualdi and Cvetkovi\'c \cite[Definition 8.1.1 and Theorem 1.2.3]{BruCve}, a nonnegative matrix
 \ $\bB = (b_{i,j})_{i,j=1,\ldots p}$ \ is irreducible provided
 that its digraph (directed graph) \ $D(\bB)$ \ (having \ $p$ \ vertices
 labeled by the numbers \ $1,2,\ldots,p$ \ and an edge from vertex \ $i$ \ to vertex \ $j$ \ provided \ $b_{i,j} > 0$)
 \ is strongly connected (that is, for each pair \ $i$ \ and
 \ $j$ \ of distinct vertices, there is a path from \ $i$ \ to \ $j$ \ and a
 path from \ $j$ \ to \ $i$).
\ Now \ $\alpha_p > 0$ \ implies that \ $D(\bA)$ \ contains a cycle
 \ $1 \to p \to (p-1) \to \cdots \to 2 \to 1$, \ hence \ $D(\bA)$ \ is strongly
 connected.

\noindent Using that \ $\bA$ \ is nonnegative and irreducible, by Horn and Johnson \cite[Theorem 8.4.4]{HJ},
 we have \ $\varrho(\bA)$ \ is an eigenvalue of \ $\bA$ \ and hence
 \begin{align*}
   \varrho(\bA)^{p} - \alpha_1 \varrho(\bA)^{p-1}-\cdots -\alpha_{p-1} \varrho(\bA) - \alpha_p = 0,
 \end{align*}
 which yields \eqref{char_poli0}.
Since
 \[
  \varphi^\prime(\lambda) = p\lambda^{p-1} - (p-1)\alpha_1\lambda^{p-2} - \cdots - \alpha_{p-1},
    \qquad \lambda\in\CCb,
 \]
 we have
 \begin{align*}
  \varphi^\prime(\varrho(\bA))
    & = p\varrho(\bA)^{-1} \sum_{k=1}^p \alpha_k\varrho(\bA)^{p-k}
       - \sum_{k=1}^{p-1} (p-k)\alpha_k\varrho(\bA)^{p-k-1} \\
   & = \sum_{k=1}^p k \alpha_k \varrho(\bA)^{p-k-1}
     = \varrho(\bA)^{p-1} \sum_{k=1}^p k \alpha_k \varrho(\bA)^{-k},
 \end{align*}
 which yields \eqref{char_poli1}.

Further, \eqref{char_poli0} yields that
 \begin{align*}
      \text{if}\qquad
          & \varrho(\bA)\,
             \begin{cases}
              < & \\
              = & \\
              >
             \end{cases}
              1,
              \qquad \text{then} \qquad
              1 = \sum_{k=1}^p \alpha_k \varr(\bA)^{-k}
              \begin{cases}
               > & \\
                = & \\
               <
              \end{cases}
               \sum_{k=1}^p \alpha_k.
        \end{align*}
This readily implies \eqref{char_poli2}.

\noindent {(ii):}
By Brualdi and Cvetkovi\'c \cite[Definition 8.2.1 and Theorem 8.2.7]{BruCve}, an irreducible
 nonnegative matrix \ $\bB = (b_{i,j})_{i,j=1,\ldots p}$ \ is primitive provided that the index of
 imprimitivity of \ $\bB$ \ (the greatest common divisor of the lengths of
 the cycles of its digraph \ $D(\bB)$) \ equals \ $1$.
\ Now the cycles of \ $D(\bA)$ \ are \ $1 \to i \to (i-1) \to \cdots \to 2 \to 1$
 \ for all \ $i=1,\ldots,p$ \ such that \ $\alpha_i > 0$ \ (not considering rotations).
 Since such a cycle has length \ $i$, \ we get the index of imprimitivity of \ $\bA$
 \ is \ $d=1$, \ which yields that \ $\bA$ \ is primitive.

\noindent The other assertions of (ii) except the uniqueness of \ $\bu_\bA$ \ and \ $\bv_\bA$
 \ follows by the Frobenius-Perron theorem, see, e.g., Horn and Johnson \cite[Theorems 8.2.11 and 8.5.1]{HJ}.
The uniqueness of \ $\bu_\bA$ \ follows by Horn and Johnson \cite[Corollary 8.2.6]{HJ} using that
 \ $\varrho(\bA^m)=\varrho(\bA)^m$ \ for all \ $m\in\NN$.
\ The uniqueness of \ $\bv_\bA$ \ can be checked as follows.
Using that the irreducibility and primitivity of \ $\bA$ \ yields the irreducibility and primitivity
 of \ $\bA^\top$ \ (see, e.g., page 507 in Horn and Johnson \cite{HJ}),
 by Horn and Johnson \cite[Theorems 8.2.11, 8.5.1 and Corollary 8.2.6]{HJ} we get
 \ $\varr(\bA^\top)=\varr(\bA)$ \ is an eigenvalue of \ $\bA^\top$, \ the algebraic and
 geometric multiplicity of \ $\varr(\bA)$ \ equal 1, corresponding to the eigenvalue
 \ $\varr(\bA)$ \ there exists a unique vector \ $\widetilde\bv_\bA \in \RR^p$ \ with positive
 coordinates such that \ $\bA^\top \widetilde\bv_\bA = \varrho(\bA) \widetilde\bv_\bA$ \ and the sum of
 the coordinates of \ $\widetilde\bv_A$ \ is 1.
Further, by Horn and Johnson \cite[page 501, Problem 1]{HJ}, we also have \ $\bu_\bA^\top \widetilde\bv_\bA>0$.
\ Using that the geometric multiplicity of \ $\varr(\bA^\top)=\varr(\bA)$ \ equals 1,
 we get \ $\bv_\bA:=\frac{1}{\bu_\bA^\top \widetilde\bv_\bA}\widetilde\bv_\bA$ \ is a unique vector
 with positive coordinates such that  $\bA^\top \bv_\bA = \varrho(\bA) \bv_\bA$ \ and
  \ $\bu_\bA^\top \bv_\bA=1$.

The forms of \ $\bu_\bA$ \ and \ $\bv_\bA$ \ can be checked as follows.
Using that they are unique it remains to verify that the imposed conditions are satisfied
 by the given forms.
We easily have \ $\bu_\bA$ \ has positive coordinates of which the sum is 1.
Further, with the notation \ $\bA\bu_\bA = \big[ (\bA\bu_\bA)_1,\ldots,(\bA\bu_\bA)_p \big]^\top$,
 \ we get
 \begin{align*}
   (\bA\bu_\bA)_1 & = \sum_{i=1}^p \alpha_i \bu_{\bA,i}
                    = \frac{\sum_{i=1}^p \alpha_i\varrho(\bA)^{-i+1}}
                     {\sum_{k=1}^p \varrho(\bA)^{-k+1}}
                    = \frac{\varrho(\bA)}{\sum_{k=1}^p \varrho(\bA)^{-k+1}}
                      \sum_{i=1}^p \alpha_i\varrho(\bA)^{-i}
                    = \frac{\varrho(\bA)}{\sum_{k=1}^p \varrho(\bA)^{-k+1}} \\
                  & = \varrho(\bA) \bu_{\bA,1},
 \end{align*}
 where the last but one equality follows by \eqref{char_poli0}.
Similarly, for \ $i=2,\ldots,p$, \ we get
 \begin{align*}
    (\bA\bu_\bA)_i = \bu_{\bA,i-1}
                   = \frac{\varrho(\bA)^{-i+2}}{\sum_{k=1}^p \varrho(\bA)^{-k+1}}
                   = \varrho(\bA) \bu_{\bA,i}.
 \end{align*}

\noindent Moreover, we easily have \ $\bv_\bA$ \ has positive coordinates and
 \begin{align*}
   \bu_\bA^\top\bv_\bA
    & = \frac{1}{\varrho(\bA)^{-p+1}\varphi^\prime(\varrho(\bA))}
       \sum_{i=1}^p\left(\varrho(\bA)^{-i+1}\sum_{\ell=i}^p\alpha_\ell\varrho(\bA)^{i-1-\ell}\right) \\
    & = \frac{1}{\varrho(\bA)^{-p+1}\varphi^\prime(\varrho(\bA))}
       \sum_{i=1}^p\sum_{\ell=i}^p \alpha_\ell\varrho(\bA)^{-\ell}
     =1,
 \end{align*}
 where the last equality follows by \eqref{char_poli1}.
With the notation
 $\bA^\top\bv_\bA = \big[ (\bA^\top\bv_\bA)_1,\ldots,(\bA^\top\bv_\bA)_p \big]^\top$, we get
 for all \ $i=1,\ldots,p-1$,
 \begin{align*}
    (\bA^\top\bv_\bA)_i & = \alpha_i \bv_{\bA,1} + \bv_{\bA,i+1}
                          = \frac{\sum_{k=1}^p
                          \varrho(\bA)^{-k+1}}{\varrho(\bA)^{-p+1}\varphi^\prime(\varrho(\bA))}
                             \left(\alpha_i\sum_{\ell=1}^p \alpha_\ell \varrho(\bA)^{-\ell}
                                + \sum_{\ell=i+1}^p \alpha_\ell \varrho(\bA)^{i-\ell}
                               \right)\\
                        & = \frac{\sum_{k=1}^p
                        \varrho(\bA)^{-k+1}}{\varrho(\bA)^{-p+1}\varphi^\prime(\varrho(\bA))}
                          \left(\alpha_i + \sum_{\ell=i+1}^p \alpha_\ell \varrho(\bA)^{i-\ell} \right)
                        = \frac{\sum_{k=1}^p
                        \varrho(\bA)^{-k+1}}{\varrho(\bA)^{-p+1}\varphi^\prime(\varrho(\bA))}
                           \sum_{\ell=i}^p \alpha_\ell \varrho(\bA)^{i-\ell}\\
                       & = \varrho(\bA) \bv_{\bA,i}.
 \end{align*}
Finally, using that  \ $\sum_{k=1}^p \alpha_k\varrho(\bA)^{-k} = 1$, \ we get
  \begin{align*}
     (\bA^\top\bv_\bA)_p
       = \alpha_p \bv_{\bA,1}
       = \alpha_p \frac{\sum_{k=1}^p \varrho(\bA)^{-k+1}}{\varrho(\bA)^{-p+1}\varphi^\prime(\varrho(\bA))}
       = \frac{\sum_{k=1}^p \varrho(\bA)^{-k+1}}{\varrho(\bA)^{-p+1}\varphi^\prime(\varrho(\bA))}
          \varrho(\bA)\alpha_p \varrho(\bA)^{-1}
       = \varrho(\bA)\bv_{\bA,p}.
  \end{align*}
\proofend

\begin{Rem}
          If \ $\alpha_p > 0$, \ $d=1$ \ and \ $\varr(\bA)=1$, \ then the unique vectors
           \ $\bu_\bA$ and \ $\bv_\bA$ \ defined in (ii) of Proposition \ref{Proposition2}
           take the forms
             \ $ \bu_\bA = \frac{1}{p} \bone_p$ \ with
             \ $\bone_p := \begin{bmatrix} 1  , \ldots , 1 \\ \end{bmatrix}^\top \in \RR^{p \times 1},$
             \ and
           \begin{align*}
             \bv_\bA = \frac{p}{\alpha_1 + 2 \alpha_2 + \cdots + p \alpha_p}
                       \begin{bmatrix}
                        \alpha_1 + \alpha_2 + \cdots + \alpha_p \\
                        \alpha_2 + \cdots + \alpha_p \\
                        \vdots \\
                        \alpha_p
                       \end{bmatrix} .
           \end{align*}
\proofend
\end{Rem}

\begin{Def}\label{DEF_primitive}
An \INARp \ process \ $(X_n)_{n\geq -p+1}$ \ with coefficients \ $\alpha_1,\ldots,\alpha_p$
 \ is called primitive if
  \renewcommand{\labelenumi}{{\rm(\roman{enumi})}}
 \begin{enumerate}
   \item $\alpha_p>0$,
   \item $d=1$, \ where \ $d$ \ is the greatest common divisor of the set
         \ $\big\{i\in\{1,\ldots,p\} : \alpha_i>0 \big\}$.
 \end{enumerate}
\end{Def}

\begin{Rem}\label{REMARK3}
If \ $\alpha_p = 0$ \ and there exists \ $i\in\{1,\ldots,p\}$ \ such that \ $\alpha_i>0$,
 \ then \ $(X_n)_{n\geq -p+1}$ \ is an INAR($p^\prime$) process with coefficients
 \ $\alpha_1,\ldots,\alpha_{p^\prime}$ \ with \ $\alpha_{p^\prime}>0$, \ where
 \ $p^\prime=\max\{i\in\{1,\ldots,p\} : \alpha_i>0\}$.
If \ $\alpha_p>0$, \ but \ $d\geq 2$, \ then the process takes the form
 \[
    X_k = \alpha_d \circ X_{k-d} + \cdots
          + \alpha_{(p/d-1)d} \circ X_{k - (p/d-1)d} + \alpha_{p} \circ X_{k - p} + \vare_k ,
    \qquad k \in \NN ,
 \]
 and hence the subsequences \ $(X_{dn-j})_{n\geq -p/d+1}$, \ $j=0,1,\ldots,d-1$, \ form independent
 primitive INAR($p/d$) processes with coefficients \ $\alpha_d,\alpha_{2d},\ldots,\alpha_p$ \ such that
 \ $X_{-p+d-j}=X_{-p+2d-j}=\cdots=X_{-j}=0$.
\ Note also that in this case not all of the coefficients \ $\alpha_d,\alpha_{2d},\ldots,\alpha_p$
 \ are necessarily positive.
Finally, we remark that an \INARp\ process \ $(X_n)_{n\geq -p+1}$ \ is primitive if and only if
 its matrix \ $\bA$ \ defined in \eqref{HELP_MATRIX_IRASMOD} is primitive.
Indeed, if \ $(X_n)_{n\geq -p+1}$ \ is primitive, then part (ii) of Proposition \ref{Proposition2}
 readily yields that \ $\bA$ \ is primitive.
Conversely (using the notations of the proof of Proposition \ref{Proposition2}), if \ $\bA$ \ is primitive,
 then, by the proof of part (i) of Proposition \ref{Proposition2}, the digraph \ $D(\bA)$ \ is strongly
 connected.
This yields that \ $\alpha_p>0$, \ since otherwise there would be no path from \ $1$ \ to \ $p$.
\ Further, the primitivity of \ $\bA$ \ yields that the index of imprimitivity of \ $\bA$ \ equals $1$.
\ Using that the cycles of \ $D(\bA)$ \ are \ $1 \to i \to (i-1) \to \cdots \to 2 \to 1$
 \ for all \ $i=1,\ldots,p$ \ such that \ $\alpha_i > 0$ \ (not considering rotations)
 and such a cycle has length \ $i$, \ we get \ $d=1$.
\proofend
\end{Rem}

The next proposition is about the limit behavior of \ $\EE(X_k)$ \ as \ $k\to\infty$.
\ This proposition can also be considered as a motivation for the classification of
 \INARp \ processes, see later on.

\begin{Pro}\label{Proposition1}
Let \ $(X_n)_{n\geq -p+1}$ \ be an \INARp \ process such that \ $X_0=X_{-1}=\cdots=X_{-p+1}=0$
 \ and \ $\EE(\vare_1^2)<\infty$.
\ Then the following assertions hold:
 \renewcommand{\labelenumi}{{\rm(\roman{enumi})}}
 \begin{enumerate}
   \item If \ $\varrho(\bA) < 1$, \ then
          \begin{align*}
            \lim_{k\to\infty} \EE(X_k)
              = \frac{\mu_\vare}{1-\sum_{i=1}^p\alpha_i}.
           \end{align*}
   \item If \ $\varrho(\bA) = 1$, \ then
         \[
            \lim_{k\to\infty} k^{-1}\EE(X_k)
               = \frac{\mu_\vare}{\sum_{i=1}^pi\alpha_i}
               =\frac{\mu_\vare}{\varphi^\prime(1)},
         \]
         where \ $\varphi$ \ is the characteristic polynomial of the matrix \ $\bA$ \ defined in
         \eqref{HELP_MATRIX_IRASMOD}.
   \item If \ $\varrho(\bA) > 1$, \ then
         \[
            \lim_{k\to\infty} \varrho(\bA)^{-kd}\EE(X_{kd-j})
              = \frac{d \mu_\vare}
                 {(\varrho(\bA)^d-1)\sum_{k=1}^p k\alpha_k\varrho(\bA)^{-k}}
              = \frac{d \mu_\vare \varrho(\bA)^{p-1}}{(\varrho(\bA)^d-1) \varphi^\prime(\varrho(\bA))}
         \]
         for all \ $j=0,1,\ldots,d-1$, \ where \ $d$ \ is the greatest common divisor of the set
         \ $\big\{i\in\{1,\ldots,p\} : \alpha_i>0 \big\}$.
  \end{enumerate}
\end{Pro}

\noindent{\bf Proof.}
If \ $\alpha_1=\cdots=\alpha_p=0$, \ then \ $\varrho(\bA)=0$ \ and \ $X_k=\vare_k$, $k\in\NN$,
 \ which yields that \ $\lim_{k\to\infty}\EE(X_k) = \mu_\vare$, \ i.e., part (i) is satisfied
  in the case of \ $\alpha_1=\cdots=\alpha_p=0$.
\ If not all of the coefficients  \ $\alpha_1,\ldots,\alpha_p$ \ are 0, then, by Remark \ref{REMARK3},
 \ $(X_n)_{n\geq -p+1}$ \ is an INAR($p^\prime$) process where
 \ $p^\prime=\max\big\{ i\in\{1,\ldots,p\} : \alpha_i>0\big\}$.
\ Hence in what follows we may and do suppose that the original process  \ $(X_n)_{n\geq -p+1}$
 \ is such that \ $\alpha_p>0$.

First we prove the proposition in the case of \ $\alpha_p>0$ \ and \ $d=1$, \ i.e.,
 in the case of \ $(X_n)_{n\geq -p+1}$ \ is primitive.

\noindent\textsf{ Proof of (i) in the case of \ $\alpha_p>0$ \ and \ $d=1$:}
In this case we verify that
 \begin{align*}
 \lim_{k\to\infty} \EE(X_k)
              = \mu_\vare \be_1^\top \sum_{j=0}^\infty \bA^j \be_1
              = \mu_\vare \be_1^\top (\bI_{\mathbf p} - \bA)^{-1} \be_1
              = \frac{\mu_\vare}{1-\alpha_1-\cdots-\alpha_p}.
 \end{align*}
By \eqref{expect_rec}, it is enough to prove that
 if \ $\varrho(\bA)<1$, \ then the series \ $\sum_{j=0}^\infty \bA^j$ \ is convergent
  and its sum is \ $(\bI_{\mathbf p} - \bA)^{-1}$.
\ By \eqref{rate}, we have
 \begin{align*}
  \sum_{j=0}^\infty \Vert \bA^j\Vert
    \leq \sum_{j=0}^\infty \varrho(\bA)^j
      \left(\Vert \varrho(\bA)^{-j} \bA^j - \bPi_\bA\Vert + \Vert \bPi_\bA\Vert\right)
    \leq  \sum_{j=0}^\infty \varrho(\bA)^j  c_\bA r_\bA^j
           + \sum_{j=0}^\infty \varrho(\bA)^j \Vert \bPi_\bA\Vert
    <\infty,
 \end{align*}
 since \ $\varrho(\bA)<1$ \ and \ $r_\bA<1$.
\ One can give another proof for the convergence of  \ $\sum_{j=0}^\infty \Vert\bA^j\Vert$.
\ Indeed, by Horn and Johnson \cite[Corollary 5.6.14]{HJ}, we have
 \ $\varrho(\bA)=\lim_{n\to\infty}\Vert A^n\Vert^{1/n}$ \ and hence
 comparison test yields the assertion.
Finally, by Lemma 5.6.10 and Corollary 5.6.16 in Horn and Johnson \cite{HJ},
 we have \ $\sum_{j=0}^\infty \bA^j = (\bI_{\mathbf p} - \bA)^{-1}$, \ and hence,
 by Cramer's rule,
 \begin{align*}
   \be_1^\top (\bI_{\mathbf p} - \bA)^{-1} \be_1
     = \frac{1}{\det(\bI_{\mathbf p} - \bA)}
     = \frac{1}{\varphi(1)}
     = \frac{1}{1-\alpha_1-\cdots-\alpha_p}.
 \end{align*}

\noindent\textsf{Proof of (ii) in the case of \ $\alpha_p>0$ \ and \ $d=1$:}
In this case we verify that
         \[
            \lim_{k\to\infty} k^{-1}\EE(X_k) = \mu_\vare \be_1^\top\bPi_\bA\,\be_1
               = \frac{\mu_\vare}{\sum_{i=1}^pi\alpha_i }
               =\frac{\mu_\vare}{\varphi^\prime(1)}.
         \]
By \eqref{expect_rec}, we get
\begin{align*}
    \EE(X_k)
      & = \mu_\vare \be_1^\top \sum_{j=0}^{k-1} \bA^j \be_1
        = \mu_\vare \be_1^\top  \sum_{j=0}^{k-1}
        \big(\bPi_\bA + (\bA^j - \bPi_\bA)\big)\be_1 \\
      & = k \mu_\vare \be_1^\top \bPi_\bA \be_1
        + \mu_\vare \be_1^\top  \sum_{j=0}^{k-1}
              (\bA^j - \bPi_\bA)\be_1,
        \qquad k\in\NN.
 \end{align*}
By \eqref{rate}, we have
 \[
    \sum_{j=0}^\infty \Vert\bA^j - \bPi_\bA\Vert
    \leq \sum_{j=0}^\infty c_\bA r_\bA^j<\infty,
 \]
 which yields that
 \[
   \lim_{k\to\infty}\frac{1}{k}\sum_{j=0}^{k-1}(\bA^j - \bPi_\bA) =  \mathbf 0,
 \]
 where \ $\mathbf 0$ \ denotes the \ $p\times p$ \ nullmatrix.
This implies \ $\lim_{k\to\infty} k^{-1}\EE(X_k) = \mu_\vare \be_1^\top\bPi_\bA\,\be_1$.
By Proposition \ref{Proposition2},  in the case of \ $\alpha_p>0$ \ and \ $d=1$ \
 ($\varrho(\bA)$ \ is not necessarily 1) we have
  \begin{align}\label{HELP_EX_CONVER}
    \be_1^\top\bPi_\bA\,\be_1
       = \be_1^\top \bu_\bA\bv_\bA^\top\,\be_1
       = \bu_{\bA,1} \bv_{\bA,1}
       =  \frac{\sum_{\ell=1}^p \alpha_\ell\varrho(\bA)^{-\ell} }
               {\varrho(\bA)^{-p+1}\varphi^\prime(\varrho(\bA))}
       = \frac{\varrho(\bA)^{p-1}}{\varphi^\prime(\varrho(\bA))}.
  \end{align}
By \eqref{char_poli2}, we have \ $\alpha_1+\cdots+\alpha_p=1$, \ and hence
 \begin{align*}
   \be_1^\top\bPi_\bA\,\be_1
     & = \frac{1}{\varphi^\prime(1)}
      = \frac{1}{p-(p-1)\alpha_1-(p-2)\alpha_2-\cdots-2\alpha_{p-2}-\alpha_{p-1}} \\
     & = \frac{1}{\sum_{i=1}^{p-1}i\alpha_i + p\left(1- \sum_{i=1}^{p-1}\alpha_i\right)},
 \end{align*}
 which yields part (ii) in the case of \ $\alpha_p>0$ \ and \ $d=1$.

\noindent\textsf{Proof of (iii) in the case of \ $\alpha_p>0$ \ and \ $d=1$:}
In this case we verify that
 \begin{align*}
  \lim_{k\to\infty} \varrho(\bA)^{-k}\EE(X_k)
             & = \frac{\mu_\vare}{\varrho(\bA)-1}\be_1^\top\bPi_\bA\,\be_1
               = \frac{\mu_\vare}{(\varrho(\bA)-1)\sum_{k=1}^p k\alpha_k\varrho(\bA)^{-k}} \\
             & = \frac{\mu_\vare \varrho(\bA)^{p-1}}{(\varrho(\bA)-1) \varphi^\prime(\varrho(\bA))}.
  \end{align*}
By \eqref{expect_rec}, we get for all \ $k\in\NN$,
\begin{align*}
   \varrho(\bA)^{-k}\EE(X_k)
     & =  \varrho(\bA)^{-k} \mu_\vare \be_1^\top \sum_{j=0}^{k-1} \bA^j \be_1
       =  \varrho(\bA)^{-k} \mu_\vare \be_1^\top  \sum_{j=0}^{k-1}\big(\varrho(\bA)^{j}\bPi_\bA
          + (\bA^j - \varrho(\bA)^{j}\bPi_\bA)\big)\be_1 \\
     & = \mu_\vare \be_1^\top \sum_{j=0}^{k-1} \varrho(\bA)^{j-k} \bPi_\bA\be_1
         + \mu_\vare \be_1^\top \varrho(\bA)^{-k} \sum_{j=0}^{k-1}(\bA^j - \varrho(\bA)^j\bPi_\bA)\be_1.
 \end{align*}
Since \ $\varrho(\bA)^{-1}<1$, \ we have
 \[
    \sum_{j=0}^{k-1} \varrho(\bA)^{j-k}
      =  \sum_{\ell=1}^{k} (\varrho(\bA)^{-1})^{\ell}
      \to \frac{1}{\varrho(\bA)-1}
      \qquad \text{as \ $k\to\infty$.}
 \]
Further, by \eqref{rate}, for all \ $k\in\NN$,
 \begin{align*}
   \left\Vert\varrho(\bA)^{-k} \sum_{j=0}^{k-1}(\bA^j - \varrho(\bA)^j\bPi_\bA)\right\Vert
        \leq \sum_{j=0}^{k-1} \varrho(\bA)^{-k+j} \Vert \varrho(\bA)^{-j}\bA^j - \bPi_\bA\Vert
        \leq c_\bA \sum_{j=0}^{k-1} \varrho(\bA)^{-k+j} r_\bA^j.
 \end{align*}
If \ $\varrho(\bA)r_\bA\ne 1$, \ then
 \[
    \left\Vert\varrho(\bA)^{-k} \sum_{j=0}^{k-1}(\bA^j - \varrho(\bA)^j\bPi_\bA)\right\Vert
       \leq c_\bA\frac{\varrho(\bA)^{-k}-r_\bA^k}{1-\varrho(\bA)r_\bA}
       \to 0 \qquad  \text{as \ $k\to\infty$,}
 \]
 since \ $\varrho(\bA)>1$ \ and \ $r_\bA<1$.
\ If \ $\varrho(\bA)r_\bA = 1$, \ then
 \[
   \left\Vert\varrho(\bA)^{-k} \sum_{j=0}^{k-1}(\bA^j - \varrho(\bA)^j\bPi_\bA)\right\Vert
      \leq c_\bA\frac{k}{\varrho(\bA)^k}
         \to 0  \qquad  \text{as \ $k\to\infty$.}
 \]
Using also \eqref{char_poli1} and \eqref{HELP_EX_CONVER}, this concludes (iii) in the case
 of \ $\alpha_p>0$ \ and \ $d=1$.

Now we turn to give a proof in the case of \ $\alpha_p>0$ \ and \ $d\geq 2$.
\ In this case, by Proposition \ref{Proposition2}, \ $\bA$ \ is irreducible, \ $\varrho(\bA)>0$ \ and,
 by Remark \ref{REMARK3}, the subsequences \ $(X_{dn-j})_{n\geq -p/d+1}$, \ $j=0,1,\ldots,d-1$, \ form
 independent primitive INAR($p/d$) processes with coefficients \ $\alpha_d,\alpha_{2d},\ldots,\alpha_p$
 \ such that \ $X_{-p+d-j}=X_{-p+2d-j}=\cdots=X_{-j}=0$.
\ Let us introduce the matrix
 \[
      \btA
        := \begin{bmatrix}
              \alpha_d & \alpha_{2d} & \alpha_{3d} & \cdots & \alpha_{p-d} & \alpha_p \\
                  1    &     0   &     0    & \cdots &      0      &     0    \\
                  0    &     1   &     0    & \cdots &      0      &     0    \\
              \vdots   & \vdots  & \vdots  &  \ddots & \vdots      & \vdots   \\
                  0    &     0   &     0    & \cdots &      0      &     0    \\
                  0    &     0   &     0    & \cdots &      1      &     0
          \end{bmatrix}
         \in\RR_+^{(p/d)\times (p/d)},
 \]
 and its characteristic polynomial
 \[
    \widetilde\varphi(\lambda)
       := \det(\lambda \bI_{\mathbf {p/d}}-\btA)
        = \lambda^{p/d} - \alpha_d \lambda^{p/d-1} - \alpha_{2d} \lambda^{p/d-2}
          - \cdots - \alpha_{p-d}\lambda -\alpha_p,
        \qquad \lambda\in\CCb.
 \]
Since the greatest common divisor of the set
 \ $\big\{ i\in\{1,\ldots,p/d\} : \alpha_{id}>0\big\}$ \ is 1, by Proposition \ref{Proposition2},
 we have \ $\btA$ \ is primitive.
We check that \ $\varrho(\bA)^d = \varrho(\btA)$.
\ Since \ $\varphi(\lambda) = \widetilde\varphi(\lambda^d)$, $\lambda\in\CCb$, \ we get
 \ $\varrho(\bA)^d \leq \varrho(\btA)$.
\ By Proposition \ref{Proposition2}, \ $\varrho(\btA)>0$ \ and \ $\varrho(\btA)$ \ is an eigenvalue
 of \ $\btA$.
\ Hence \ $\varrho(\btA)^{1/d}$ \ is an eigenvalue of \ $\bA$, \ which implies that
 \ $\varrho(\bA)\geq \varrho(\btA)^{1/d}$ \ or \ equivalently \ $\varrho(\bA)^d\geq \varrho(\btA)$.

If \ $\varrho(\bA)<1$, \ then \ $\varrho(\btA)<1$ \ and using that part (i) has already been proved for
 primitive matrices (i.e., in the case of \ $\alpha_p>0$ \ and \ $d=1$) \ we have for all
 \ $j=0,1,\ldots,d-1$,
 \[
   \lim_{n\to\infty} \EE(X_{nd-j})
     = \frac{\mu_\vare}{1-\alpha_d-\alpha_{2d}-\cdots - \alpha_p}
     = \frac{\mu_\vare}{1- \sum_{i=1}^p\alpha_i}.
 \]
This yields that \ $\lim_{n\to\infty}\EE(X_n)$ \ exists with the given limit in (i).

If \ $\varrho(\bA)=1$, \ then \ $\varrho(\btA)=1$ \ and using that part (ii) has already been proved for
 primitive matrices we have for all \ $j=0,1,\ldots,d-1$,
 \begin{align*}
   \lim_{n\to\infty}\frac{\EE(X_{dn-j})}{n}
      = \frac{\mu_\vare}{\alpha_d+2\alpha_{2d}+\cdots +\frac{p}{d}\alpha_p}
      = \frac{d\mu_\vare}{d\alpha_d+2d\alpha_{2d}+\cdots + p\alpha_p}
      = \frac{d\mu_\vare}{\sum_{i=1}^pi\alpha_i}.
 \end{align*}
This yields that \ $\lim_{k\to\infty}k^{-1}\EE(X_k)$ \ exists with given limit in (ii).

If \ $\varrho(\bA)>1$, \ then \ $\varrho(\btA)>1$ \ and using that part (iii) has already been proved for
 primitive matrices we have for all \ $j=0,1,\ldots,d-1$,
 \begin{align*}
  \lim_{n\to\infty}\frac{\EE(X_{nd-j})}{\varrho(\btA)^n}
     & = \frac{\mu_\vare}{(\varrho(\btA)-1)\sum_{k=1}^{p/d} k\alpha_{kd}\varrho(\btA)^{-k}}
       = \frac{d\mu_\vare}{(\varrho(\bA)^d-1)\sum_{k=1}^{p/d} kd\alpha_{kd}\varrho(\bA)^{-kd}} \\
     &  = \frac{d\mu_\vare}{(\varrho(\bA)^d-1)\sum_{\ell=1}^p \ell \alpha_{\ell}\varrho(\bA)^{-\ell}}
        = \frac{d\mu_\vare \varrho(\bA)^{p-1}}{(\varrho(\bA)^d-1)\varphi^\prime(\varrho(\bA))},
 \end{align*}
 where the last equality follows by \eqref{char_poli1}.
Since
 \begin{align*}
   \lim_{k\to\infty}\frac{\EE(X_{kd-j})}{\varrho(\bA)^{kd}}
      = \lim_{k\to\infty}\frac{\EE(X_{kd-j})}{\varrho(\btA)^{k}},
      \qquad j=0,1,\ldots,d-1,
 \end{align*}
 we have (iii).
\proofend

Based on the asymptotic behavior of \ $\EE(X_k)$ \ as \ $k\to\infty$ \ described in Proposition \ref{Proposition1},
 we distinguish three cases.
The case \ $\varrho(\bA) < 1$ \ is called \emph{stable} or
 \emph{asymptotically stationary},
 whereas the cases \ $\varrho(\bA) = 1$ \ and \ $\varrho(\bA) > 1$ \ are called
 \emph{unstable} and \emph{explosive}, respectively.
Note also that, if \ $\alpha_p>0$, \ then, by \eqref{char_poli2} of Proposition \ref{Proposition2},
 \ $\varrho(\bA) < 1$, \ $\varrho(\bA) = 1$ \ and \ $\varrho(\bA) > 1$ \ are equivalent
 with \ $\alpha_1 + \cdots + \alpha_p < 1$, \ $\alpha_1 + \cdots + \alpha_p = 1$ \ and
 \ $\alpha_1 + \cdots + \alpha_p > 1$, \ respectively.

\section{Convergence of unstable \INARp\ processes}\label{conv_section}

A function \ $f : \RR_+ \to \RR$ \ is called \emph{c\`adl\`ag} if it is right
 continuous with left limits.
\ Let \ $\DD(\RR_+, \RR)$ \ and \ $\CC(\RR_+, \RR)$ \ denote the space of all
 real-valued c\`adl\`ag and continuous functions on \ $\RR_+$, \ respectively.
Let \ $\cD_\infty$ \ denote the Borel $\sigma$-field in \ $\DD(\RR_+, \RR)$
 \ for the metric defined in (16.4) in Billingsley \cite{Bil}
 (with this metric \ $\DD(\RR_+, \RR)$ \ is a complete and separable metric space).
For stochastic processes \ $(\cY_t)_{t\in\RR_+}$ \ and \ $(\cY^n_t)_{t\in\RR_+}$,
 \ $n\in\NN$, \ with c\`adl\`ag paths we write \ $\cY^n \distr \cY$ \ if the
 distribution of \ $\cY^n$ \ on the space \ $(\DD(\RR_+, \RR),\cD_\infty)$ \ converges weakly to
 the distribution of \ $\cY$ \ on the space \ $(\DD(\RR_+, \RR),\cD_\infty)$ \ as \ $n\to\infty$.

For each \ $n \in \NN$, \ consider the random step processes
 \[
   \cX^n_t := n^{-1} X_{\nt} , \qquad t \in \RR_+ , \quad n \in \NN ,
 \]
 where \ $\lfloor x\rfloor$ \ denotes the integer part of a real number \ $x\in\RR$.
The positive part of \ $x \in \RR$ \ will be denoted by \ $x^+$.

\begin{Thm}\label{main}
Let \ $(X_k)_{k \geq -p+1}$ \ be a primitive \INARp\ process with coefficients
 \ $\alpha_1,\ldots,\alpha_p \in [0,1]$ \ such that
 \ $\alpha_1 + \cdots + \alpha_p = 1$ \ (hence it is unstable).
Suppose that \ $X_0=X_{-1}=\cdots=X_{-p+1}=0$ \ and \ $\EE(\vare_1^2) < \infty$.
\ Then
 \begin{equation}\label{Conv_X}
  \cX^n \distr \cX \qquad \text{as \ $n \to \infty$,}
 \end{equation}
 where \ $(\cX_t)_{t\in\RR_+}$ \ is the unique strong solution of the stochastic differential
 equation (SDE)
 \begin{equation}\label{CIR2}
   \DS \dd \cX_t
   = \frac{1}{\varphi^\prime(1)}
     \Big( \mu_\vare \, \dd t
           + \sqrt{\sigma_\alpha^2 \cX_t^+} \,\dd \cW_t \Big) ,
  \qquad t\in\RR_+,
 \end{equation}
 with initial value \ $\cX_0 = 0$, \ where
 \[
   \varphi^\prime(1) = \alpha_1 + 2 \alpha_2 + \cdots + p \alpha_p > 0, \qquad
   \sigma_\alpha^2 := \alpha_1(1-\alpha_1) + \cdots + \alpha_p(1-\alpha_p) ,
 \]
 and \ $(\cW_t)_{t\in\RR_+}$ \ is a standard Wiener process.
 (Here \ $\varphi$ \ is the characteristic polynomial of the matrix \ $\bA$ \ defined in
 \eqref{HELP_MATRIX_IRASMOD}.)
\end{Thm}

\begin{Rem}
Note that under the conditions Theorem \ref{main}, if \ $p\geq 2$, then \ $\sigma_\alpha^2>0$,
 \ and if \ $p=1$, \ then \ $\sigma_\alpha^2=0$.
Indeed, if \ $p\geq 2$, \ then \ $\alpha_p<1$, \ since otherwise \ $\alpha_1=\cdots=\alpha_{p-1}=0$
 \ and hence the greatest common divisor of \ $\{i\in\{1,\ldots,p\} : \alpha_i>0\} = \{ p\}$ \
 would be \ $p$, \ which is a contradiction.
Since, by our assumption \ $\alpha_p>0$, \ we get \ $\sigma_\alpha^2\geq \alpha_p(1-\alpha_p)>0$. \
If \ $p=1$, \ then \ $\alpha_p=\alpha_1=1$, \ and hence \ $\sigma_\alpha^2=\alpha_1(1-\alpha_1)=0$.

\noindent Remark also that in the case of \ $p=1$ \ we have \ $\alpha_1=1$ \ and hence
 \ $X_n=\sum_{i=1}^n\vare_i$, $n\in\NN$, \ $\varphi^\prime(1)=1$, \ $\sigma_\alpha^2=0$ \ and then
 the limit process in Theorem \ref{main}
 is deterministic, namely \ $\cX_t=\mu_\vare t$, $t\in\RR_+$.
\ To describe the asymptotic behavior of an unstable INAR$(1)$ process one has to go one step further and
 one has to investigate the fluctuation limit.
By Donsker's theorem (see, e.g., Billingsley \cite[Theorem 8.2]{Bil}), we have
 \ $\sqrt{n}(\cX^n-\EE(\cX^n))\distr \sigma_\vare \cW$ \ as \ $n\to\infty$, \ where \ $\cW$ \ is
 a standard Wiener process.
For completeness, we remark that Isp\'any, Pap and Zuijlen \cite[Proposition 4.1]{IspPapZui1}
 describes the fluctuation limit behavior of nearly unstable INAR($1$) processes.
 \proofend
\end{Rem}

\begin{Rem}\label{REMARK1}
The SDE \eqref{CIR2} has a unique strong solution \ $(\cX_t^x)_{t\geq 0}$ \ for all
 initial values  \ $\cX_0^x=x\in\RR$.
Indeed, since \ $\vert \sqrt{x} - \sqrt{y}\vert\leq \sqrt{\vert x-y\vert}$, $x,y\geq 0$,
 \ the coefficient functions
 \ $\RR\ni x \mapsto \mu_\vare/\varphi^\prime(1)$ \ and
 \ $\RR\ni x \mapsto \sqrt{\sigma_\alpha^2 x^+}/\varphi^\prime(1)$ \ satisfy
 conditions of part (ii) of Theorem 3.5 in Chapter IX in Revuz and Yor \cite{RevYor}
 or the conditions of Proposition 5.2.13 in Karatzas and Shreve \cite{KarShr}.
Further, by the comparison theorem (see, e.g., Revuz and Yor \cite[Theorem 3.7, Chapter IX]{RevYor}),
 if the initial value \ $\cX_0^x=x$ \ is nonnegative, then \ $\cX_t^x$ \ is nonnegative
 for all \ $t\in\RR_+$ \ with probability one.
Hence \ $\cX_t^+$ \ may be replaced by \ $\cX_t$ \ under the square root in \eqref{CIR2}.
The unique strong solution of the SDE \eqref{CIR2} is known as a squared Bessel process,
 a squared-root process or a Cox-Ingersoll-Ross (CIR) process.
 \proofend
\end{Rem}

\begin{Rem}
If the matrix \ $\bA$ \ is not primitive but unstable, then we can suppose that
 \ $\alpha_p>0$, \ since otherwise it is an unstable INAR($p^\prime$) process with
 \ $p^\prime:=\max\big\{i\in\{1,\ldots,p\} : \alpha_i>0 \big\}$ \ (note that there exists
 an \ $i\in\{1,\ldots,p\}$ \ such that \ $\alpha_i>0$ \ because of the unstability of \ $\bA$).
\ If \ $\alpha_p>0$ \ and \ $d\geq 2$, \ then, by Remark \ref{REMARK3}, the subsequences
 \ $(X_{dn-j})_{n\geq -p/d+1}$, \ $j=0,1,\ldots,d-1$, \ form independent primitive unstable INAR($p/d$)
 processes with coefficients \ $\alpha_d,\alpha_{2d},\ldots,\alpha_p$ \ such that
 \ $X_{-p+d-j}=X_{-p+2d-j}=\cdots=X_{-j}=0$.
\ Hence one can use Theorem \ref{main} for these subsequences.
With the notations
 \[
   \cX^{n,j}_t := \frac{1}{n} X_{d\lfloor nt\rfloor - j},\qquad t\in\RR_+,\;\;
                      n\geq - \frac{p}{d} +1,\;\;j=0,1,\ldots,d-1,
 \]
 by Theorem \ref{main}, \ $\cX^{n,j}\distr \cX^{(j)}$ \ as \ $n\to\infty$, \ where
  \ $(\cX^{(j)}_t)_{t\in\RR_+}$ \ is the unique strong solution of the SDE
 \begin{align*}
   \DS \dd \cX^{(j)}_t
   = \frac{1}{\alpha_d+2\alpha_{2d}+\cdots+ \frac{p}{d} \alpha_p}
     \Big( \mu_\vare \, \dd t
           + \sqrt{\sigma_\alpha^2 (\cX^{(j)}_t)^+} \,\dd \cW^{(j)}_t \Big) ,
  \qquad t\in\RR_+,
 \end{align*}
 with initial value \ $\cX^{(j)}_0 = 0$ \ and \ $(\cW^{(j)}_t)_{t\in\RR_+}$, \ $j=0,1,\ldots,d-1$,
 \ are independent standard Wiener processes.
We note that if \ $\alpha_p>0$ \ and \ $d\geq 2$, \ then
 \ $\cX^n$ \ does not converge in general as \ $n \to \infty$.
By giving a counterexample, we show that even the 2-dimensional distributions do not converge in general.
Let \ $p:=4$, \ $\alpha_1=\alpha_3:=0$, \ $\alpha_2=\alpha_4:=1/2$.
\ Then \ $d=2$ \ and using that \ $\cX^{n,j}\distr \cX^{(j)}$ \ as \ $n\to\infty$, \ $j=0,1$,
 \ we have
 \begin{align}\label{HELP_COUNTEREXAMPLE}
   [\cX^{n,\,0}_1,\cX^{n,\,0}_2] = \left[\frac{1}{n}X_{2n}, \frac{1}{n}X_{4n} \right]
     \quad\text{converges in distribution to}\quad
   [\cX^{(0)}_1,\cX^{(0)}_2]
 \end{align}
 as \ $n\to\infty$, \ and
 \begin{align}\label{HELP_COUNTEREXAMPLE2}
   [\cX^{n,1}_1,\cX^{n,1}_2] = \left[\frac{1}{n}X_{2n-1}, \frac{1}{n}X_{4n-1} \right]
     \quad\text{converges in distribution to}\quad
   [\cX^{(1)}_1,\cX^{(1)}_2]
 \end{align}
 as \ $n\to\infty$, \ where \ $(\cX^{(j)}_t)_{t\in\RR_+}$ \ is the unique strong solution of the SDE
  \begin{align*}
   \DS \dd \cX^{(j)}_t
   = \frac{2}{3}
     \Big( \mu_\vare \, \dd t
           + \sqrt{\frac{1}{2}(\cX^{(j)}_t)^+} \,\dd \cW^{(j)}_t \Big) ,
  \qquad t\in\RR_+,
 \end{align*}
 with initial value \ $\cX^{(j)}_0 = 0$, \ $j=0,1$.
\ However, we show that
 \begin{align}\label{HELP_COUNTEREXAMPLE3}
   [\cX^{n}_1,\cX^{n}_2]
     = \left[\frac{1}{n}X_{n}, \frac{1}{n}X_{2n} \right]
       \quad \text{does not converge in distribution as \ $n\to\infty$.}
 \end{align}
Indeed, we have
 \begin{align*}
   [\cX^{2n}_1,\cX^{2n}_2]
     = \left[\frac{1}{2n}X_{2n}, \frac{1}{2n}X_{4n} \right]
      =  \left[\frac{1}{2}\cX^{n,\,0}_1,\frac{1}{2} \cX^{n,\,0}_2\right]
 \end{align*}
 and hence, by \eqref{HELP_COUNTEREXAMPLE},
 \begin{align*}
   [\cX^{2n}_1,\cX^{2n}_2]
       \quad \text{converges in distribution to} \quad
       \left[\frac{1}{2}\cX^{(0)}_1,\frac{1}{2}\cX^{(0)}_2\right]
       \qquad \text{as \ $n\to\infty$.}
 \end{align*}
Further, using that
 \[
   [\cX^{2n-1}_1,\cX^{2n-1}_2]
      = \left[\frac{1}{2n-1}X_{2n-1}, \frac{1}{2n-1}X_{2(2n-1)} \right]
      =  \left[\frac{n}{2n-1}\cX^{n,1}_1,\cX^{2n-1,\,0}_1\right]
 \]
 and that the subsequences \ $(X_{2n-1})_{n\in\NN}$ \ and \ $(X_{2(2n-1)})_{n\in\NN}$ \ are
 independent, by \eqref{HELP_COUNTEREXAMPLE} and \eqref{HELP_COUNTEREXAMPLE2}, we get
 \begin{align*}
   [\cX^{2n-1}_1,\cX^{2n-1}_2]
     \quad \text{converges in distribution to} \quad
   \left[\frac{1}{2}\cX^{(1)}_1, \cX^{(0)}_1\right]
    \qquad \text{as \ $n\to\infty$.}
 \end{align*}
Since the random variables
 \[
     \left[\frac{1}{2}\cX^{(0)}_1,\frac{1}{2}\cX^{(0)}_2\right]
     \qquad \text{and} \qquad
     \left[\frac{1}{2}\cX^{(1)}_1, \cX^{(0)}_1\right]
 \]
 do not have the same distributions (the coordinates of the first one are dependent,
 however the coordinates of the second one are independent), we get \eqref{HELP_COUNTEREXAMPLE3}.
\proofend
\end{Rem}

For proving Theorem \ref{main}, let us introduce the sequence
 \begin{equation}\label{Mk}
  M_k := X_k - \EE(X_k \mid \cF_{k-1})
       = X_k - \alpha_1 X_{k-1} - \cdots - \alpha_p X_{k-p} - \mu_\vare ,
  \qquad k \in \NN ,
 \end{equation}
 of martingale differences with respect to the filtration
 \ $(\cF_k)_{k \in \ZZ_+}$, \ and the random step processes
 \[
   \cM^n_t := n^{-1} \sum_{k=1}^{\nt} M_k ,
   \qquad t \in \RR_+ , \quad n \in \NN .
 \]
First we will verify convergence
 \begin{equation}\label{Conv_M}
   \cM^n \distr \cM \qquad \text{as \ $n \to \infty$,}
 \end{equation}
 where \ $(\cM_t)_{t\in\RR_+}$ \ is the unique strong solution of the SDE
 \begin{equation}\label{CIR1}
   \DS \dd \cM_t
   = \sqrt{\frac{\sigma_\alpha^2}{\varphi^\prime(1)}(\cM_t + \mu_\vare t)^+}
     \,\dd \cW_t,
   \qquad t\in\RR_+,
 \end{equation}
 with initial value \ $\cM_0 = 0$.
\ The proof of \eqref{Conv_M} can be found in Section \ref{proof_section}.

\begin{Rem}\label{REMARK2}
 \ If \ $(\cX_t^{x})_{t\in\RR_+}$ \ is a strong solution of \eqref{CIR2} with initial value
 \ $\cX_0^{x}=x\in\RR$, \ then, by It\^o's formula, \ $\cM_t^{x} := \varphi^\prime(1) \cX_t^{x} - \mu_\vare t$,
 \ $t \in \RR_+$, \ is a strong solution of \eqref{CIR1} with initial value \ $\cM_0^{x} = \varphi'(1)x$. \
On the other hand, if \ $(\cM_t^{y})_{t\in\RR_+}$ \ is a strong solution of
 \eqref{CIR1} with initial value \ $\cM_0^{y}=y\in\RR$, \ then, again by It\^o's formula,
 \begin{align}\label{HELP_ITO}
   \cX_t^{y} := \frac{1}{\varphi^\prime(1)}(\cM_t^y + \mu_\vare t),\qquad t \in \RR_+,
 \end{align}
 is a strong solution of \eqref{CIR2} with initial value \ $\cX_0^{y}=\frac{1}{\varphi'(1)}y$.
\ Hence, by Remark \ref{REMARK1}, the SDE \eqref{CIR1} has a unique strong solution
 \ $(\cM_t^{y})_{t\geq 0}$ \ for all initial values \ $\cM_0^{y}=y\in\RR$.
\ Further, if the initial value \ $\cM_0^{y}=y$ \ is nonnegative, then
 \ $\cM_t^{y}+\mu_\vare t$ \ is nonnegative for all \ $t\in\RR_+$ \ with probability one.
Hence \ $(\cM_t + \mu_\vare t)^+$ \ may be replaced by \ $(\cM_t + \mu_\vare t)$
 \ under the square root in \eqref{CIR1}.
 \proofend
\end{Rem}

Moreover, from \eqref{Mk} we obtain the recursion
 \begin{align}\label{recursion}
  X_k = \alpha_1 X_{k-1} + \cdots + \alpha_p X_{k-p} + M_k + \mu_\vare ,
  \qquad k \in \NN ,
 \end{align}
 which can be written in the form
 \ $\bX_k = \bA \bX_{k-1} + (M_k + \mu_\vare) \be_1$, \ $k \in \NN$.
\ Consequently,
 \[
   \bX_k = \sum_{j=1}^k (M_j + \mu_\vare) \bA^{k-j} \be_1 ,
   \qquad k \in \NN ,
 \]
 implying
 \begin{equation}\label{REGR}
   X_k = \be_1^\top \bX_k
       = \sum_{j=1}^k (M_j + \mu_\vare) \be_1^\top \bA^{k-j} \be_1 ,
   \qquad k \in \NN .
 \end{equation}
In Section \ref{proof_section}, we show that the statement \eqref{Conv_X} will follow
 from \eqref{Conv_M} and \eqref{REGR} using a version of the continuous mapping theorem
 (see Appendix).

\section{Application to Boston armed robberies data set}\label{Boston_armrob}

\begin{figure}[thpb]
      \centering
      \includegraphics[height=4cm, width=8cm]{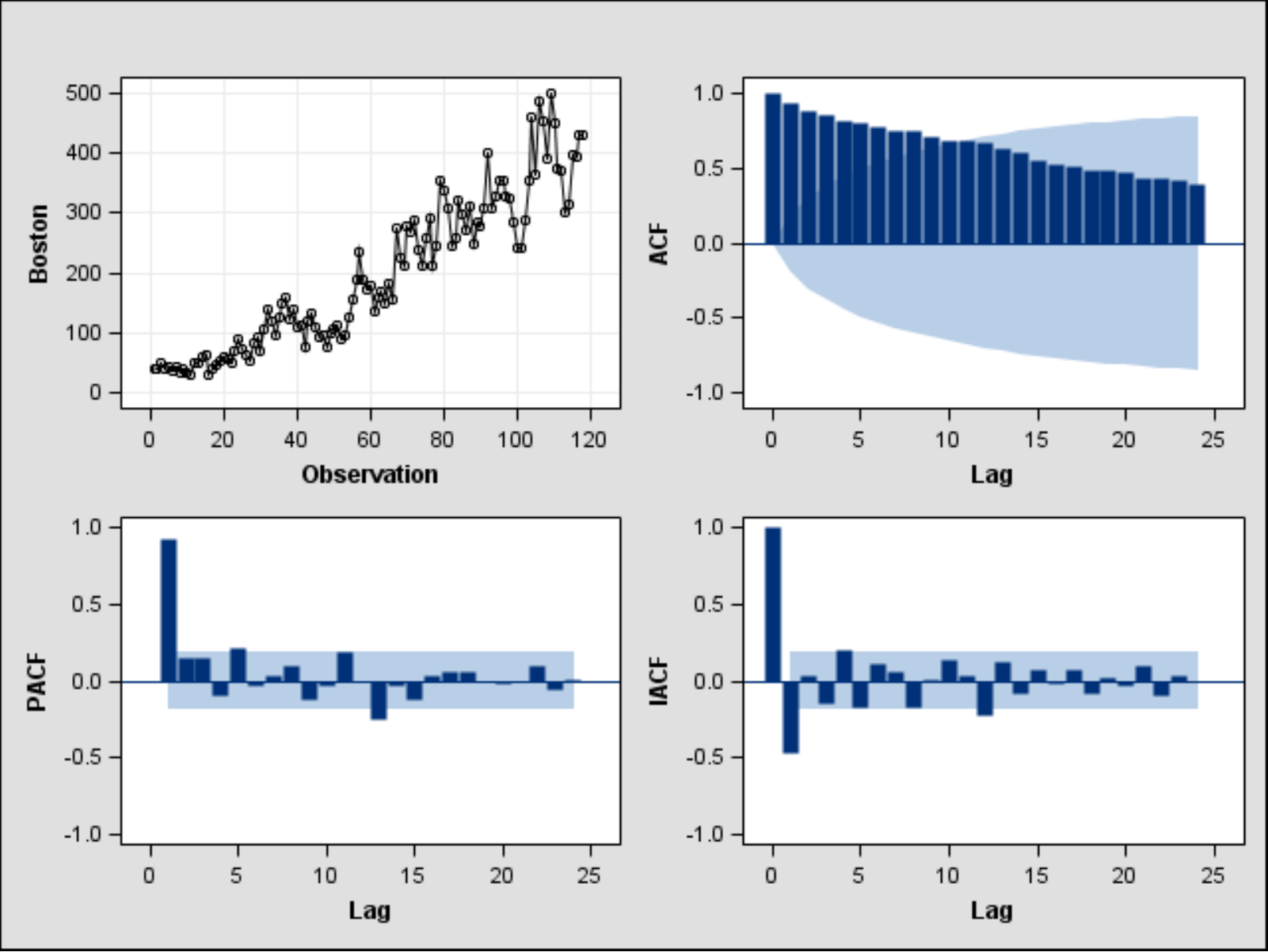}
      \caption{Boston armed robberies, time series (top left),
       autocorrelation function (top right), partial autocorrelation function
       (bottom left), inverse autocorrelation function (bottom right).}
      \label{boston}
\end{figure}

This data set consists of 118 counts of monthly armed robberies in Boston from
 January 1966 to October 1975 (Fig.~\ref{boston}). The data were originally
 published in Deutsch and Alt \cite{DA}, see also the time series 6.10 in O'Donovan
 \cite[Appendix A.3]{OD}.
It can also be obtained from the Time Series Data Library:
 http://robjhyndman.com/tsdldata/data/mccleary5.dat.
Deutsch and Alt \cite{DA} used this time series to illustrate the method of intervention analysis
 developed by Box and Tiao \cite{BT}.
They assessed the impact of a 1975 Massachusetts gun control law on armed robbery in Boston.
 The correlation analysis for this series, shown in Fig.~\ref{boston}, and preliminary ARIMA model
 fitting clearly indicate unstability.
For preliminary fitting, subset ARIMA$(12,0,0)$ (Model 1) and ARIMA$(1,0,0)
\times (1,0,0)_{12}$ (Model 2) models are applied. Here and in the sequel,
ARIMA$(p,d,q)\times(P,D,Q)_s$ denotes a seasonal ARIMA model with period
$s\in\NN$ and orders $(p,d,q),(P,D,Q)\in\ZZ_+^3$, where capital letters
denote the seasonal orders. We use the following approach to characterize
the unstability of an ARIMA model. Let $a(B)$ and $A(B)$ be the
autoregressive and seasonal autoregressive polynomial of the model, respectively,
 where $B$ denotes the backshift operator. We suppose that these
polynomials are stable, i.e., the roots are all lie outside the complex unit circle.
Define the coefficients $\alpha_i$, $i=1,\ldots, p+d+s(P+D)$,
by $a(B)A(B)(1-B)^d(1-B^s)^D=1-\sum_{i=1}^{p+d+s(P+D)}
\alpha_i B^i$. Then, we characterize the unstability of the model by the sum
$\Sigma:=\sum_{i=1}^{p+d+s(P+D)}\alpha_i$. Clearly, if an ARIMA model is unstable
 (nonstationary), i.e., $d>0$ or $D>0$, and hence its characteristic
 polynomial has unit root 1, then $\Sigma=1$.
Since Model 1 is unstable and Model 2 is nearly unstable, see Table \ref{model_table}, Deutsch and Alt
\cite{DA} suggested first order differencing and seasonal differencing getting an
ARIMA$(0,1,1)\times(0,1,1)_{12}$ model (Model 3).
\begin{table}[h]
{\small
   \begin{tabular}{c c c c}
      \hline\hline
      Model & Fitted model & $\Sigma$ & Standard error \\
      \hline\hline
      1 & $(1-0.7865 B-0.2135 B^{12}) X_k =  \vare_k +116.3733$ & 1 & 39.55 \\
      2 & $(1-0.9783B)(1-0.2677B^{12}) X_k = \vare_k +49.2087$ & 0.9841 & 40.39 \\
      3 & $(1-B)(1-B^{12})X_k = (1-0.5154B)(1-0.7345B^{12})
            \vare_k +0.3181 $ & 1&38.66 \\
      4 & $(1-B)\ln X_k = (1-0.4345B)(1+0.1886B^{12}) \vare_k +0.0195$ & 1 &
      0.1954 \\
      5 &$ X_k=0.6069\circ X_{k-1}+0.412\circ X_{k-12} +14.971+ \tvare_k$ &1.0189& 526.8\\
      6 &$X_k=0.682\circ X_{k-1}+0.3497\circ X_{k-12} +9.961+ \tvare_k$ &1.0317& 26.18\\
      \hline
   \end{tabular}}
\caption{Fitted models for Boston armed robberies data set with $\Sigma$ and
standard error.}
\label{model_table}
\end{table}
In contrast, Hay and McCleary
\cite{HM} claimed that Deutsch and Alt had misspecified the stochastic component for
this time series and they proposed only first order differencing getting an
ARIMA$(0,1,1)\times (0,0,1)_{12}$ model (Model 4) after logarithmic transformation
of the time series. Hay and McCleary reported that this alternative model
has better statistical properties and there is no intervention into the
time series (i.e., the parameters of the model do not vary in time),
 thus there is inconclusive evidence for the effect claimed by Deutsch and Alt.
They argued for the need of logarithmic transformation to eliminate the ``variance''
 nonstationarity of the time series.
The following was reported by Hay and McCleary \cite{HM}:
``We conducted several analyses to obtain supporting evidence for our hypothesis of
variance nonstationarity. First, we divided the series into equal length segments
and calculated the mean and standard deviation for each segment. Both statistics
showed a nearly monotonic increase over time and were highly intercorrelated. Two
tests of homogeneity of variance (Cochran's C and the Bartlett and Box's F) also
indicated that the segment variances were not homogeneous.''

Based on the foregoing it is evident that the Boston armed robberies data set
possesses the following properties: it is integer--valued, heteroscedastic,
and unstable. Our aim here is to fit an appropriate INAR(p) model for this data set
using the method of conditional least squares (CLS) and to compare our model with
the previously mentioned ones. The CLS estimators $\halpha_i$, $i=1,\ldots,p$,
and $\hmue$ of the parameters $\alpha_i$, $i=1,\ldots,p$, and $\mu_\vare$ of an INAR$(p)$
model based on the observations $X_1,\ldots,X_n$ are given by minimizing the residual
sum of squares $\sum_{k=p+1}^n M_k^2$ in \eqref{recursion}. This technique has
been suggested by Klimko and Nelson \cite{KN} for general stochastic processes,
and it has been applied for INAR$(p)$ models by Du and Li \cite[Theorem 4.2]{DuLi}
 proving the asymptotic normality of these estimators in the stable case.
The correlation analysis (Fig.~\ref{boston}) shows that there are significant
 dependences between $X_k$ and $X_{k-1}$, and, due to the seasonal effect, between
$X_k$ and $X_{k-12}$. Thus, we fit a subset INAR$(12)$ model where the strictly
positive coefficients are $\alpha_1$ and $\alpha_{12}$, and we estimate these
  (autoregressive) parameters and the mean $\mu_\vare$. By solving the normal equations
we have Model 5, see Table \ref{model_table}, where $\tvare_k:= \vare_k - \hmue$
 is the centered innovation.
Similarly to ARIMA models we characterize the unstability of an INAR$(p)$ model
 by the sum $\Sigma:= \sum_{i=1}^{p}\alpha_i$ (the classification of INAR(p) models is based
 on this sum, see the end of Section \ref{INARp_section}).
Then the fitted Model 5 appears to be unstable since $\Sigma=1.0189$. For the goodness--of--fit of ARIMA
 and INAR models the standard error (the square root of the mean square error) is applied which is
defined by $\SE:=((n-p-r)^{-1}\sum_{k=p+1}^n \hM_k^2)^{1/2}$, where
$\hM_k:=X_k-\sum_{i=1}^p \halpha_i X_{k-i} - \hmue$, $k=p+1,\ldots,n$, are the estimated
residuals and $r$ denotes the number of estimated parameters. The standard error is
relatively high for Model 5 ($\SE= 526.8$) comparing with that of  Deutsch and Alt's
model (Model 3) because the ``error'' terms $M_k$ fluctuate to much in \eqref{recursion}
if the INAR model is unstable.
(We note that the model of Hay and McCleary (Model 4) is uncomparable with the other ones
using the standard error because of the logarithmic transformation has changed the scale.)

\begin{figure}[thpb]
      \centering
      \includegraphics[height=4cm, width=8cm]{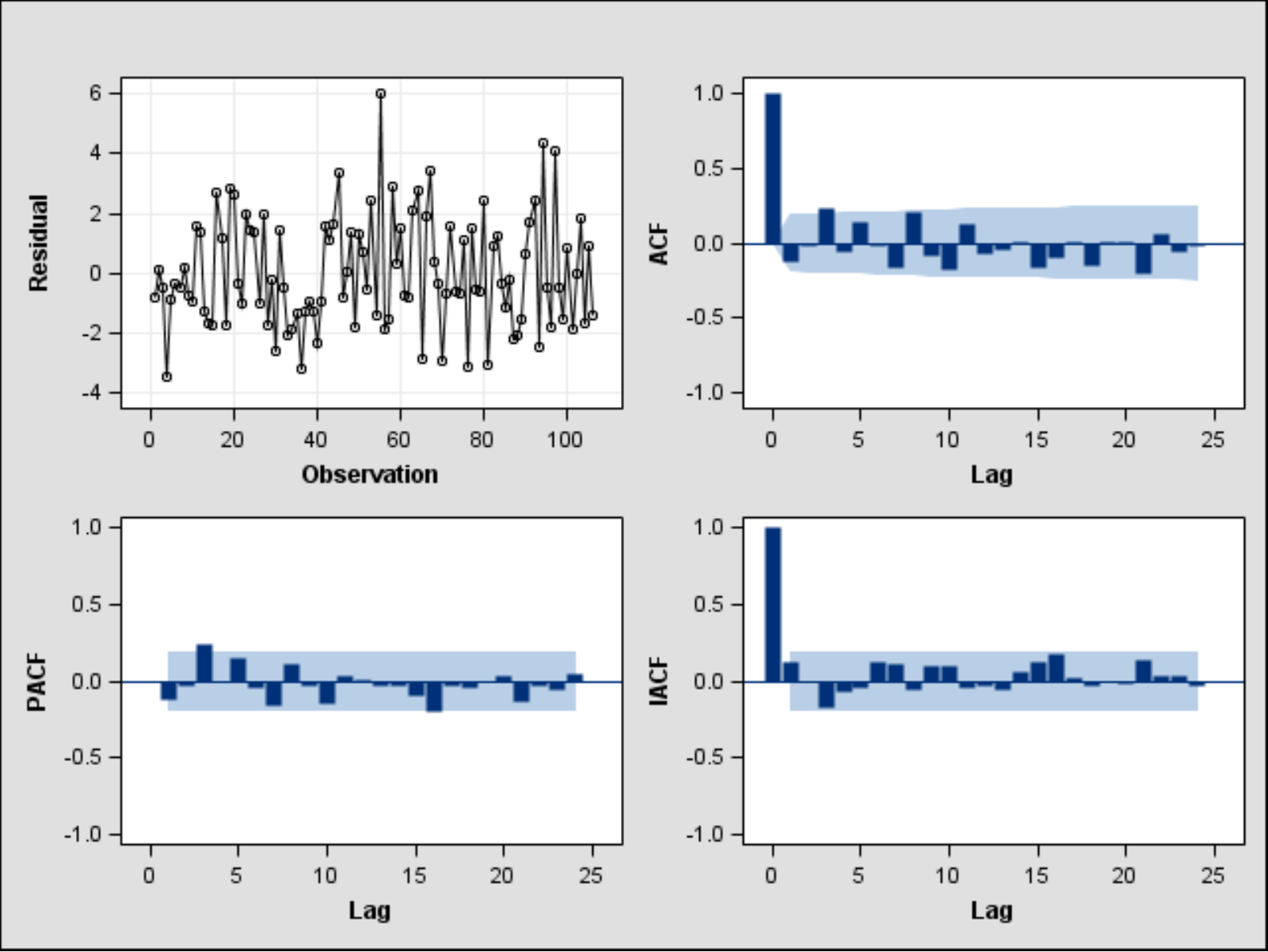}
      \caption{Residual analysis of Model 6, residual series (top left),
       autocorrelation function (top right), partial autocorrelation function
       (bottom left), inverse autocorrelation function (bottom right).}
      \label{residual}
\end{figure}

To stabilize the fluctuation of $M_k$ let us introduce the weighted martingale differences
 \[
   M_k^w:=\frac{M_k}{\left(\sum_{\{j:\alpha_j>0\}} X_{k-j}+1\right)^{1/2}}, \qquad k=p+1,\ldots,n.
 \]
Note that $\EE(M_k^w\,|\,\cF_{k-1})=0$ and, by \eqref{Mcond},
$$
     \EE\left((M_k^w)^2\bmid \cF_{k-1}\right)= \frac{\sum_{\{j:\alpha_j>0\}}\alpha_j(1-\alpha_j)
      X_{k-j}+\sigma_\vare^2}{\sum_{\{j:\alpha_j>0\}} X_{k-j}+1},
      \qquad k=p+1,\ldots,n.
$$
Since $\EE\left((M_k^w)^2 \bmid \cF_{k-1}\right)\le \sum_{\{j:\alpha_j>0\}}\alpha_j(1-\alpha_j)
+\sigma_\vare^2$, the conditional variance of the ``weighted error'' terms $M_k^w$ would
not fluctuate too much even if $(X_k)_{k\in\NN}$ is unbounded.
Moreover, we have $\EE\left((M_k^w)^2
 \bmid \cF_{k-1}\right)\to \frac{1}{c}\sum_j \alpha_j(1-\alpha_j)$ almost surely as $X_k\to\infty$
 and $X_k/X_{k-1}\to 1$ almost surely, where \ $c$ \ denotes the cardinality of the set
 \ $\{j\in\{1,\ldots,p\} : \alpha_j>0\}$. \ Hence, the weighted error terms $M_k^w$ are
 asymptotically homogeneous in the stable and the unstable cases as well.
The weighted conditional least squares
(WCLS) estimation is given by minimizing the weighted residual sum of squares
$\sum_{k=p+1}^n (M_k^w)^2$. This technique has been suggested by Wei and Winnicki
\cite{WW2} for branching processes with immigration to derive a unified estimation
procedure for the offspring mean. By solving the normal equations we have Model 6 which
appears to be unstable again, see Table \ref{model_table}. Defining the standard error
for Model 6 as $\SE:=((n-p-r)^{-1}\sum_{k=p+1}^n (\hM_k^w)^2)^{1/2}$, this
subset INAR(12) model possesses the smallest standard error among the fitted models except
that of Hay and McCleary. The correlation analysis of estimated weighted residuals
$\hM_k^w$, see Fig.~\ref{residual}, shows that they form a white noise time series.

In summary, Model 6 is an adequate model for Boston armed robberies times series
 since its coefficients can be considered significant, it has minimum number of parameters
 and minimal residual variance (among the fitted models), and the residuals form a white noise.
We note that the asymptotic theory of CLS and WCLS estimation of INAR$(p)$ models in the
unstable case has not yet been developed now, this is a task for the future.
Finally, we would like to call attention to other possible estimation methods which may also
 work in the unstable case.
For example, Enciso--Mora et al.~\cite{ENR} proposed a reversible jump MCMC algorithm which
even works well near the borders of the stationary region and has been successfully
applied to a simulated nearly unstable INAR(3) data set having \ $\Sigma=0.99$ \
 as the sum of the (autoregressive) parameters.


\section{Proof of Theorem \ref{main}}\label{proof_section}

For the proof we will use Corollary \ref{EEX}, Theorem \ref{Conv2DiffCor}
 and Lemma \ref{Conv2Funct} which can be found in Appendix.

First we prove \eqref{Conv_M}, i.e., \ $\cM^n\distr\cM$ \ as \ $n\to\infty$.
\ We will apply Theorem \ref{Conv2DiffCor} for \ $\cU = \cM$, \ $U^n_k = n^{-1} M_k$, \ $n,\,k \in \NN$,
 \ and for \ $(\cF_k^n)_{k\in\ZZ_+}=(\cF_k)_{k\in\ZZ_+}$, $n\in\NN$.
\ By Remark \ref{REMARK2}, the SDE \eqref{CIR1} has
 a unique strong solution for all initial values \ $\cM^{x}_0=x$, \ $x\in\RR$.
\ Now we show that conditions (i) and (ii) of Theorem \ref{Conv2DiffCor} hold.
We have to check that for each \ $T>0$,
 \begin{align} \label{Cond1}
  &\sup_{t\in[0,T]}
    \bigg| \frac{1}{n^2} \sum_{k=1}^{\nt}
            \EE( M_k^2 \mid \cF_{k-1} )
            - \frac{\sigma_\alpha^2}{\varphi^\prime(1)}
              \int_0^t (\cM^n_s + \mu_\vare s)^+ \, \dd s \, \bigg|
   \stoch 0,\\
  &\frac{1}{n^2}
   \sum_{k=1}^{\nT}
    \EE( M_k^2 \bbone_{\{|M_k| > n\theta\}} \mid \cF_{k-1} \big)
   \stoch 0
   \qquad\text{for all \ $\theta>0$} \label{Cond2}
 \end{align}
 as \ $n\to\infty$, \ where \ $\stoch$ \ means convergence in probability.

By \eqref{Mk} and using also that \ $\alpha_1+\cdots+\alpha_p=1$, \ we get
 \begin{align*}
  \cM^n_s + \mu_\vare s
  & = n^{-1}
      \sum_{k=1}^{\ns}
       \left( X_k - \sum_{i=1}^p\alpha_i X_{k-i} - \mu_\vare \right)
      + \mu_\vare s \\
  & = n^{-1}
      \left( \sum_{k = \ns - p +1}^{\ns} X_k
             - \sum_{i=1}^{p-1} \alpha_i \sum_{k = \ns - p +1}^{\ns - i} X_k
      \right)
      + \frac{ns - \ns}{n} \mu_\vare \\
  & = \frac{1}{n}
       \sum_{j=1}^{p} \sum_{i=j}^{p} \alpha_i X_{\ns - j +1}
      + \frac{ns - \ns}{n} \mu_\vare .
 \end{align*}

Thus \ $(\cM^n_s + \mu_\vare s)^+ = \cM^n_s + \mu_\vare s$, \ and using that
 \[
   \int_0^t \frac{ns - \ns}{n}\,\dd s
     = \frac{t^2}{2} - \frac{1}{n}\left(\frac{1}{n}\sum_{k=1}^{\nt-1}k
                                        +\left(t-\frac{\nt}{n}\right)\nt\right)
     =\frac{\nt+(nt-\nt)^2}{2n^2},
 \]
 we get
 \begin{align*}
   \int_0^t (\cM^n_s + \mu_\vare s)^+ \,\dd s
   &= \frac{1}{n^2}
      \sum_{\ell=0}^{\nt-1}
         \sum_{j=1}^{p} \sum_{i=j}^{p} \alpha_i X_{\ell - j +1}
      + \frac{nt-\nt}{n^2}
        \sum_{j=1}^{p} \sum_{i=j}^{p} \alpha_i X_{\nt - j +1}
         \\
   & \phantom{=\:}
      +\frac{\nt+(nt-\nt)^2}{2n^2} \mu_\vare.
  \end{align*}
Hence, using that \ $\varphi^\prime(1) = \alpha_1 + 2 \alpha_2 + \cdots + p \alpha_p$, \ we have
 \begin{align*}
  \int_0^t (\cM^n_s + \mu_\vare s)^+ \,\dd s
   &= \frac{\varphi^\prime(1)}{n^2} \sum_{\ell=0}^{\nt-1} X_\ell
      - \frac{1}{n^2} \sum_{i=2}^{p} \alpha_i \sum_{j=\nt -i +1}^{\nt -1} X_{j}
       \\
   & \phantom{=\:\;}
      + \frac{nt-\nt}{n^2}
        \sum_{j=1}^{p} \sum_{i=j}^{p} \alpha_i X_{\nt - j +1}
      +\frac{\nt+(nt-\nt)^2}{2n^2} \mu_\vare .
 \end{align*}
Using \eqref{Mcond}, we obtain
 \begin{align*}
   \frac{1}{n^2}
    \sum_{k=1}^{\nt}
     \EE( M_k^2 \mid \cF_{k-1})
   & = \frac{1}{n^2}
      \sum_{k=1}^{\nt}
       \left( \sum_{i=1}^p \alpha_i(1-\alpha_i) X_{k-i} + \sigma_\vare^2 \right)\\
   &=\frac{1}{n^2}\sum_{i=1}^p\left(\alpha_i(1-\alpha_i)\sum_{j=1}^{\nt-i+1}X_{j-1}\right)
            +\frac{\nt}{n^2}\sigma_\vare^2\\[2mm]
   &= \frac{\sigma_\alpha^2}{n^2}
      \sum_{k=1}^{\nt}  X_{k-1}
        - \frac{1}{n^2} \sum_{i=2}^p \left(\alpha_i(1-\alpha_i)
                       \sum_{j=\nt -i+1}^{\nt -1} X_j\right)
      + \frac{\nt}{n^2} \sigma_\vare^2 .
 \end{align*}
Hence, for all \ $n\in\NN$, \ the randomness of the difference in \eqref{Cond1}
 is via a linear combination of the random variables \ $ X_{\nt - j}$, $j=1,\ldots,p$.
\ Then, in order to show \eqref{Cond1}, it suffices to prove
 \begin{equation}\label{Cond11}
  \sup_{t\in[0,T]}
   \frac{1}{n^2}
   X_{\nt}
  \stoch 0 \qquad \text{as \ $n\to\infty$.}
 \end{equation}
By \eqref{REGR} and \eqref{C},
 \[
   X_{\nt} \leq \sum_{j=1}^{\nt} |M_j + \mu_\vare| \cdot \| \bA^{\nt-j} \|
   \leq C_\bA \biggl( \nt \cdot \mu_\vare + \sum_{j=1}^{\nt} |M_j| \biggr) .
 \]
Consequently, in order to prove \eqref{Cond11}, it suffices to show
 \[
  \frac{1}{n^2} \sum_{j=1}^{\nT} |M_j| \stoch 0 \qquad
  \text{as \ $n\to\infty$.}
 \]
In fact, one can show that \ $n^{-2} \sum_{j=1}^{\nT} \EE(|M_j|) \to 0$.
\ Indeed, Corollary \ref{EEX} yields that
 \begin{align*}
    n^{-2} \sum_{j=1}^{\nT} \EE(|M_j|)
      \leq \frac{K}{n^2} \sum_{j=1}^{\nT} \sqrt{j}
      \leq \frac{K}{n^2} \nT\sqrt{\nT} \to 0
      \quad \text{as \ $n\to\infty$,}
 \end{align*}
 with some constant \ $K\in\RR_+$.
\ Thus we obtain \eqref{Cond1}.

To prove \eqref{Cond2}, consider the decomposition
 \ $M_k = N_k + (\vare_k - \mu_\vare)$, \ where, by \eqref{HELP_M_DECOMP},
 \begin{align*}
  N_k
  := \sum_{\ell=1}^{X_{k-1}} ( \xi_{k,1,\ell} - \EE(\xi_{k,1,\ell}) ) + \cdots
     + \sum_{\ell=1}^{X_{k-p}} ( \xi_{k,p,\ell} - \EE(\xi_{k,p,\ell}) ) .
 \end{align*}
Clearly,
 \begin{gather*}
  M_k^2
  \leq 2 \big( N_k^2 + ( \vare_k - \mu_\vare )^2 \big)
  \qquad \text{and} \qquad
  \bbone_{\{ | M_k | > n \theta \}}
  \leq \bbone_{\{ | N_k | > n \theta / 2 \}}
       + \bbone_{\{ | \vare_k - \mu_\vare | > n \theta / 2 \}} ,
 \end{gather*}
 and hence \eqref{Cond2} will be proved once we show
 \begin{gather}
  \frac{1}{n^2}
   \sum_{k=1}^{\nT}
    \EE( N_k^2 \bbone_{\{ |N_k | > n \theta \}} \mid \cF_{k-1} )
   \stoch 0 \qquad \text{for all \ $\theta>0$,} \label{Cond21} \\
  \frac{1}{n^2}
   \sum_{k=1}^{\nT}
    \EE( N_k^2 \bbone_{\{ | \vare_k - \mu_\vare | > n \theta \}} \mid \cF_{k-1})
   \stoch 0 \qquad \text{for all \ $\theta>0$,} \label{Cond22} \\
  \frac{1}{n^2}
   \sum_{k=1}^{\nT}
    \EE( (\vare_k - \mu_\vare))^2 \mid \cF_{k-1})
   \stoch 0. \label{Cond23}
 \end{gather}

First we prove \eqref{Cond21}.
Using that the random variables
 \ $\{\xi_{k,i,j} : j\in\NN, \, i\in\{1,\ldots,p\}\}$ \ are independent of the
 \ $\sigma$-algebra \ $\cF_{k-1}$ \ for all \ $k\in\NN$, \ we get
 \[
   \EE(N_k^2 \bbone_{\{ |N_k| > n\theta\}}\mid\cF_{k-1})
     = F_k(X_{k-1},\ldots,X_{k-p}),
 \]
 where \ $F_k : \ZZ_+^p \to \RR$ \ is given by
 \[
   F_k(z_1,\dots,z_p)
   :=\EE((S_k(z_1,\dots,z_p)^2 \bbone_{\{ |S_k(z_1,\dots,z_p)| > n\theta\}})),
   \qquad z_1,\ldots,z_p\in\ZZ_+,
 \]
 with
 \ $S_k(z_1,\dots,z_p)
    :=\sum_{i=1}^p \sum_{\ell=1}^{z_i} (\xi_{k,i,\ell} - \EE(\xi_{k,i,\ell}))$.
\ Consider the decomposition
 \[
   F_k(z_1,\dots,z_p) = A_k(z_1,\dots,z_p) + B_k(z_1,\dots,z_p),
 \]
 where
 \begin{align*}
  A_k(z_1,\dots,z_p)
  & := \sum_{i=1}^p \sum_{\ell=1}^{z_i}
        \EE((\xi_{k,i,\ell} - \EE(\xi_{k,i,\ell}))^2
            \bbone_{\{ |S_k(z_1,\dots,z_p)| > n\theta\}}),\\
  B_k(z_1,\dots,z_p)
  & := \sum\nolimits'
        \EE((\xi_{k,i,\ell} - \EE(\xi_{k,i,\ell}))
            (\xi_{k,j,\ell'} - \EE(\xi_{k,j,\ell'}))
            \bbone_{\{ |S_k(z_1,\dots,z_p)| > n\theta\}}) ,
 \end{align*}
 where the sum \ $\sum'$ \ is taken for \ $i,j=1,\dots,p$,
 \ $\ell=1,\dots,z_i$, \ $\ell'=1,\dots,z_j$ \ with \ $(i,\ell)\not=(j,\ell')$.
\ Consider the decompositions
 \[
   S_k(z_1,\dots,z_p)
   = (\xi_{k,i,\ell} - \EE(\xi_{k,i,\ell})) + \tS^i_{k,\ell}(z_1,\dots,z_p) ,
   \qquad i=1,\ldots,p,\;\; \ell=1,\ldots,z_i ,
 \]
 where
 \[
   \tS^i_{k,\ell}(z_1,\dots,z_p)
   :=\sum\nolimits''\big(\xi_{k,j,\ell'} - \EE(\xi_{k,j,\ell'})\big) ,
 \]
 where the sum \ $\sum''$ \ is taken for \ $j=1,\dots,p$ \ and
 \ $\ell'=1,\dots,z_j$ \ with \ $(j,\ell')\not=(i,\ell)$.

Using that
 \[
    \bbone_{\{ |S_k(z_1,\dots,z_p)| > n\theta\}}
     \leq
      \bbone_{\{ |\xi_{k,i,\ell} - \EE(\xi_{k,i,\ell})| > n\theta/2\}}
       +
      \bbone_{\{ |\tS^i_{k,\ell}(z_1,\dots,z_p)| > n\theta/2\}},
 \]
 we have
 \[
    A_k(z_1,\dots,z_p) \leq A_k^{(1)}(z_1,\dots,z_p) + A_k^{(2)}(z_1,\dots,z_p),
 \]
 where
 \begin{align*}
  A_k^{(1)}(z_1,\dots,z_p)
  & := \sum_{i=1}^p \sum_{\ell=1}^{z_i}
        \EE( (\xi_{k,i,\ell} - \EE(\xi_{k,i,\ell}))^2
             \bbone_{\{ |\xi_{k,i,\ell} - \EE(\xi_{k,i,\ell})| > n \theta/2\}}),\\
  A_k^{(2)}(z_1,\dots,z_p)
  & := \sum_{i=1}^p \sum_{\ell=1}^{z_i}
        \EE( (\xi_{k,i,\ell} - \EE(\xi_{k,i,\ell}))^2
             \bbone_{\{ |\tS^i_{k,\ell}(z_1,\dots,z_p)| > n\theta/2\}}) .
 \end{align*}
In order to prove \eqref{Cond21}, it is enough to show that
 \begin{gather}\nonumber
   \frac{1}{n^2}
    \sum_{k=1}^{\nT}
     A_k^{(1)}(X_{k-1},\ldots,X_{k-p})
    \stoch 0,\qquad\qquad
   \frac{1}{n^2}
    \sum_{k=1}^{\nT}
     A_k^{(2)}(X_{k-1},\ldots,X_{k-p})
    \stoch 0,\\ \label{Bk}
   \frac{1}{n^2}
    \sum_{k=1}^{\nT}
     B_k(X_{k-1},\dots,X_{k-p})
    \stoch 0
 \end{gather}
 as \ $n\to\infty$.
\ We have
 \[
   A_k^{(1)}(z_1,\dots,z_p)
   = \sum_{i=1}^p z_i
      \EE( (\xi_{1,i,1} - \EE(\xi_{1,i,1}))^2
           \bbone_{\{ |\xi_{1,i,1} - \EE(\xi_{1,i,1})| > n\theta/2\}}) ,
    \qquad k\in\NN,
 \]
 where
 $$
   \EE( (\xi_{1,i,1} - \EE(\xi_{1,i,1}))^2
         \bbone_{\{ |\xi_{1,i,1} - \EE(\xi_{1,i,1})| > n\theta/2\}}) \to 0,
 $$
 as \ $n\to\infty$ \ for all \ $i \in \{1,\dots,p\}$ \ by the dominated convergence theorem.
Thus, by Corollary \ref{EEX}, we get with some constant \ $K\in\RR_+$,
 \begin{align*}
   \frac{1}{n^2}\sum_{k=1}^{\nT} & \EE(A_k^{(1)}(X_{k-1},\dots,X_{k-p}))\\
      &= \frac{1}{n^2} \sum_{k=1}^{\nT} \sum_{i=1}^p \EE(X_{k-i})
                        \EE( (\xi_{1,i,1} - \EE(\xi_{1,i,1}))^2
                        \bbone_{\{ |\xi_{1,i,1} - \EE(\xi_{1,i,1})| > n\theta/2\}}) \\
     & \leq \sum_{i=1}^p \left[
                   \EE( (\xi_{1,i,1} - \EE(\xi_{1,i,1}))^2
                   \bbone_{\{ |\xi_{1,i,1} - \EE(\xi_{1,i,1})| > n\theta/2\}})
                     \frac{K}{n^2} \sum_{k=i+1}^{\nT} (k-i)
                   \right] \\
     & \leq K \frac{\nT(\nT+1)}{2n^2}
              \sum_{i=1}^p  \EE( (\xi_{1,i,1} - \EE(\xi_{1,i,1}))^2
                         \bbone_{\{ |\xi_{1,i,1} - \EE(\xi_{1,i,1})| > n\theta/2\}})
      \to 0,
 \end{align*}
  which yields \ $n^{-2} \sum_{k=1}^{\nT} A_k^{(1)}(X_{k-1},\dots,X_{k-p}) \stoch 0$.

Independence of \ $\xi_{k,i,\ell} - \EE(\xi_{k,i,\ell})$ \ and
 \ $\tS^i_{k,\ell}(z_1,\dots,z_p)$ \ implies
 \[
   A_k^{(2)}(z_1,\dots,z_p)
   = \sum_{i=1}^p \sum_{\ell=1}^{z_i}
      \EE( (\xi_{k,i,\ell} - \EE(\xi_{k,i,\ell}))^2)
      \PP( |\tS^i_{k,\ell}(z_1,\dots,z_d)| > n\theta/2\big) .
 \]
Here \ $\EE( (\xi_{k,i,\ell} - \EE(\xi_{k,i,\ell}))^2) = \alpha_i(1-\alpha_i)$, $i=1,\ldots,p$,
 \ and, by Markov's inequality,
 \begin{multline*}
  \PP( |\tS^i_{k,\ell}(z_1,\dots,z_p)| > n\theta/2)
    \leq \frac{4}{n^2\theta^2}
         \EE(\tS^i_{k,\ell}(z_1,\dots,z_p)^2) \\
  = \frac{4}{n^2\theta^2} \var(\tS^i_{k,\ell}(z_1,\dots,z_p))
    = \frac{4}{n^2\theta^2}
      \sum\nolimits'' \alpha_j(1-\alpha_j)
    \leq \frac{4}{n^2\theta^2} \sum_{j=1}^p z_j \alpha_j(1-\alpha_j).
 \end{multline*}
Thus we get
 \[
   A_k^{(2)}(z_1,\dots,z_p)
   \leq \frac{4}{n^2\theta^2} \sum_{i=1}^p \sum_{j=1}^p
        z_i z_j \alpha_i(1-\alpha_i) \alpha_j(1-\alpha_j) .
 \]
Hence, by Cauchy-Schwarz's inequality and Corollary \ref{EEX}, we get with some constant \ $K\in\RR_+$,
 \begin{align*}
   \frac{1}{n^2} \sum_{k=1}^{\nT} \EE(A_k^{(2)}(X_{k-1},\ldots,X_{k-p}))
     & \leq \frac{4}{n^4\theta} \sum_{k=1}^{\nT}  \sum_{i=1}^p \sum_{j=1}^p
             \EE(X_{k-i}X_{k-j}) \alpha_i(1-\alpha_i) \alpha_j(1-\alpha_j) \\
      & \leq \frac{4K}{n^4\theta} \sum_{k=1}^{\nT} k^2
              \left(\sum_{i=1}^p \alpha_i(1-\alpha_i)\right)^2
       \to 0,
 \end{align*}
 which implies \ $n^{-2} \sum_{k=1}^{\nT} A_k^{(2)}(X_{k-1},\ldots,X_{k-p}) \stoch 0$.

By Cauchy-Schwarz's inequality,
 \[
   |B_k(z_1,\dots,z_p)|
   \leq \sqrt{B_k^{(1)}(z_1,\dots,z_p) \,
        \EE(\bbone_{\{ |S_k(z_1,\dots,z_p)| > n\theta\}})},
 \]
 where
 \begin{multline*}
  B_k^{(1)}(z_1,\dots,z_p)
    := \EE\bigg(\Big(\sum\nolimits'
                     (\xi_{k,i,\ell} - \EE(\xi_{k,i,\ell}))
                     (\xi_{k,j,\ell'} - \EE(\xi_{k,j,\ell'})) \Big)^2\bigg),
    \qquad z_1,\ldots,z_p\in\ZZ_+.
 \end{multline*}
Using the independence of \ $\xi_{k,i,\ell} - \EE(\xi_{k,i,\ell})$ \ and
 \ $\xi_{k,j,\ell'} - \EE(\xi_{k,j,\ell'})$ \ for \ $(i,\ell) \ne (j,\ell')$, \ we get
 \begin{align*}
   B_k^{(1)}(z_1,\dots,z_p)
   & = \sum\nolimits' \alpha_i(1-\alpha_i) \alpha_j(1-\alpha_j) \\
   & = \sum_{i=1}^p z_i (z_i-1) \alpha_i^2(1-\alpha_i)^2
     + \sum_{i \ne j} z_i z_j \alpha_i(1-\alpha_i) \alpha_j(1-\alpha_j) \\
   & \leq K_1(z_1+\cdots+z_p)^2,
 \end{align*}
 with some constant \ $K_1\in\RR_+$.
\ Further, by Markov's inequality,
 \[
   \EE(\bbone_{\{ |S_k(z_1,\dots,z_p)| > n\theta\}})
   \leq \frac{1}{n^2\theta^2} \sum_{j=1}^p z_j \alpha_j(1-\alpha_j)
   \leq \frac{K_2}{n^2\theta^2}(z_1+\dots+z_p),
 \]
 with some constant \ $K_2\in\RR_+$.
\ Hence
 \[
   |B_k(z_1,\dots,z_p)|
     \leq \frac{K}{n} (z_1+\dots+z_p)^{3/2},
     \qquad z_1,\ldots,z_p\in\ZZ_+,
 \]
 with some constant \ $K\in\RR_+$.
\ Using that
 \[
   (z_1+\dots+z_p)^{3/2} \leq c_p (z_1^{3/2}+\dots+z_p^{3/2}),
      \qquad z_1,\ldots,z_p\in\ZZ_+,
 \]
 with some constant \ $c_p\in\RR_+$, \ we get, in order to show
 \eqref{Bk},  it suffices to prove
 \ $n^{-3} \sum_{k=1}^{\nT} \big( X_{k-1}^{3/2} + \cdots + X_{k-p}^{3/2} \big)
    \stoch 0$.
\ In fact,
 \ $n^{-3} \sum_{k=1}^{\nT}
           \big( \EE(X_{k-1}^{3/2}) + \cdots + \EE(X_{k-p}^{3/2}) \big) \to 0$
 \ since Corollary \ref{EEX} implies
 \ $\EE(X_\ell^{3/2}) \leq \big(\EE(X_\ell^2)\big)^{3/4} = O(\ell^{3/2})$.
\ Thus we finished the proof of \eqref{Cond21}.

Now we turn to prove \eqref{Cond22}.
Using that for all \ $k\in\NN$ \ the random variables \ $\{\xi_{k,i,j},\,\vare_k : j\in\NN, \, i\in\{1,\ldots,p\}\}$
 \ are independent of the \ $\sigma$-algebra \ $\cF_{k-1}$, \ we get\
 \ $\EE(N_k^2 \bbone_{\{ |\vare_k - \mu_\vare| > n\theta\}} \mid\cF_{k-1})
    =G_k(X_{k-1},\ldots,X_{k-p})$,
 \ where \ $G_k : \ZZ_+^p \to \RR$ \ is given by
 \[
   G_k(z_1,\dots,z_p)
   :=\EE(S_k(z_1,\dots,z_p)^2 \bbone_{\{ |\vare_k - \mu_\vare| > n\theta\}}\big),
   \qquad z_1,\ldots,z_p\in\ZZ_+.
 \]
Using again the independence of
 \ $\big\{\xi_{k,i,j}, \, \vare_k : j\in\NN, \, i\in\{1,\dots,p\}\big\}$,
 \[
   G_k(z_1,\dots,z_p)
   = \PP( |\vare_k - \mu_\vare| > n\theta)
     \sum_{i=1}^p \sum_{\ell=1}^{z_i}
      \EE( (\xi_{k,i,\ell} - \EE(\xi_{k,i,\ell}))^2),
 \]
 where by Markov's inequality,
 \ $\PP( |\vare_k - \mu_\vare| > n\theta\big)
    \leq n^{-2}\theta^{-2} \EE((\vare_k -\mu_\vare)^2)
    =  n^{-2}\theta^{-2} \sigma_\vare^2$ ,
 \ and
 \ $\EE( (\xi_{k,i,\ell} - \EE(\xi_{k,i,\ell}))^2) = \alpha_i(1-\alpha_i)$.
\ Hence, in order to show \eqref{Cond22}, it suffices to prove
 \ $n^{-4} \sum_{k=1}^{\nT} X_k \stoch 0$.
\ In fact, by Corollary \ref{EEX}, \ $n^{-4} \sum_{k=1}^{\nT} \EE(X_k) \to 0$.

Now we turn to prove \eqref{Cond23}.
By independence of \ $\vare_k$ \ and \ $\cF_{k-1}$,
 \[
   \frac{1}{n^2}
    \sum_{k=1}^{\nT}
     \EE( (\vare_k - \mu_\vare)^2 \mid \cF_{k-1} )
   = \frac{1}{n^2}
     \sum_{k=1}^{\nT}
      \EE( (\vare_k - \mu_\vare)^2 )
   = \frac{\nT}{n^2} \sigma_\vare^2
   \to 0 ,
 \]
 thus we obtain \eqref{Cond23}.
Hence we get \eqref{Cond2}, and we conclude, by
 Theorem \ref{Conv2DiffCor}, convergence \ $\cM_n \distr \cM$.

Now we start to prove \eqref{Conv_X}.
By \eqref{REGR}, \ $\cX^n = \Psi_n(\cM^n)$, \ where the mapping
 \ $\Psi_n : \DD(\RR_+, \RR) \to \DD(\RR_+, \RR)$ \ is given by
 \[
   \Psi_n(f)(t)
   := \sum_{j=1}^\nt
       \left( f\left( \frac{j}{n} \right) - f\left( \frac{j-1}{n} \right)
              + \frac{\mu_\vare}{n} \right)
       \be_1^\top \bA^{\nt - j} \be_1
 \]
 for \ $f \in \DD(\RR_+, \RR)$, \ $t \in \RR_+$, \ $n\in\NN$.
\ Further, \ $\cX = \Psi(\cM)$, \ where, by \eqref{HELP_ITO}, the mapping
 \ $\Psi : \DD(\RR_+, \RR) \to \DD(\RR_+, \RR)$ \ is given by
 \[
   \Psi(f)(t) := \frac{1}{\varphi^\prime(1)} \big( f(t) + \mu_\vare t \big) ,
   \qquad f \in \DD(\RR_+, \RR), \quad t\in\RR_+ .
 \]
We check that the mappings \ $\Psi_n$, $n\in\NN$, \ and \ $\Psi$ \ are measurable.
Continuity of \ $\Psi$ \ follows from the characterization of convergence in
 \ $\DD(\RR_+, \RR)$, \ see, e.g., Ethier and Kurtz
 \cite[Proposition 3.5.3]{EK}, thus we obtain measurability of \ $\Psi$.
\ Indeed, if \ $f_n\in\DD(\RR_+, \RR)$, $n\in\NN$, \ $f\in\DD(\RR_+, \RR)$ \ and the sequence
 \ $(f_n)_{n\in\NN}$ \ converges in \ $\DD(\RR_+, \RR)$ \ to \ $f$, \ then for all \ $T>0$
 \ there exist continuous, increasing mappings \ $\lambda_n$, $n\in\NN$, \ from
 \ $[0,\infty)$ \ onto \ $[0,\infty)$ \ such that
  \[
    \lim_{n\to\infty}\sup_{t\in[0,T]}\vert\lambda_n(t)-t\vert = 0
    \qquad \text{and}\qquad
    \lim_{n\to\infty}\sup_{t\in[0,T]}\vert f_n(\lambda_n(t))- f(t)\vert = 0.
  \]
Since for all \ $t\in\RR_+$
 \begin{align*}
   \vert \Psi(f_n)(\lambda_n(t)) - \Psi(f)(t)\vert
     & = \left\vert
            \frac{1}{\varphi'(1)}(f_n(\lambda_n(t))+\mu_\vare\lambda_n(t))
             - \frac{1}{\varphi'(1)}(f(t)+\mu_\vare t)
        \right\vert\\
     & \leq \frac{1}{\varphi'(1)}\vert f_n(\lambda_n(t)) - f(t)\vert
          +\frac{\mu_\vare}{\varphi'(1)}\vert \lambda_n(t) - t \vert,
 \end{align*}
 we have for all \ $T>0$,
 \[
   \lim_{n\to\infty}\sup_{t\in[0,T]}\vert \Psi(f_n)(\lambda_n(t)) - \Psi(f)(t)\vert = 0.
 \]
In order to prove measurability of \ $\Psi_n$, \ first we localize it.
For each \ $N \in \NN$, \ consider the stopped mapping
 \ $\Psi^N_n : \DD(\RR_+, \RR) \to \DD(\RR_+, \RR)$ \ given by
 \ $\Psi^N_n(f)(t) := \Psi_n(f)(t \land N)$ \ for \ $f \in \DD(\RR_+, \RR)$,
 \ $t \in \RR_+$, \ $n, N \in \NN$.
Obviously, \ $\Psi^N_n(f) \to \Psi_n(f)$ \ in \ $\DD(\RR_+, \RR)$ \ as
 \ $N \to \infty$ \ for all \ $f \in \DD(\RR_+, \RR)$, \ since for all
 \ $T > 0$ \ and \ $N \geq T$ \ we have \ $\Psi^N_n(f)(t) := \Psi_n(f)(t)$, $t\in[0,T]$,
 \ and hence \ $\sup_{t \in [0,T]} |\Psi^N_n(f)(t) - \Psi_n(f)(t)| \to 0$ \ as
 \ $N \to \infty$.
\ Consequently, it suffices to show measurability of \ $\Psi^N_n$ \ for all
 \ $n, N \in \NN$.
\ We can write \ $\Psi^N_n = \Psi^{N,2}_n \circ \Psi^{N,1}_n$, \ where the
 mappings \ $\Psi^{N,1}_n : \DD(\RR_+, \RR) \to \RR^{nN+1}$ \ and
 \ $\Psi^{N,2}_n : \RR^{nN+1} \to \DD(\RR_+, \RR)$ \ are defined by
 \begin{align*}
   \Psi_n^{N,1}(f)
   &:= \left( f(0), f\left( \frac{1}{n} \right), f\left( \frac{2}{n} \right), \dots, f(N) \right) , \\
   \Psi_n^{N,2}(x_0,x_1,\dots,x_{nN})(t)
   &:= \sum_{j=1}^{\lfloor n (t \land N) \rfloor}
       \left( x_j - x_{j-1} + \frac{\mu_\vare}{n} \right)
       \be_1^\top \bA^{\nt - j} \be_1
 \end{align*}
 for \ $f \in \DD(\RR_+, \RR)$, \ $t \in \RR_+$,
 \ $x = (x_0,x_1,\dots,x_{nN}) \in \RR^{nN+1}$, \ $n,N\in\NN$.
\ Measurability of \ $\Psi^{N,1}_n$ \ follows from Ethier and Kurtz
 \cite[Proposition 3.7.1]{EK}.
Next we show continuity of \ $\Psi^{N,2}_n$ \ by checking
 \ $\sup_{t \in [0,T]} |\Psi^{N,2}_n(x^k) (t)- \Psi^{N,2}_n(x)(t)| \to 0$ \ as
 \ $k \to \infty$ \ for all \ $T > 0$ \ whenever \ $x^k \to x$ \ in
 \ $\RR^{nN+1}$.
 \ This convergence follows from the estimates
 \[
   \sup_{t \in [0,T]} |\Psi^{N,2}_n(x^k) (t)- \Psi^{N,2}_n(x)(t)|
   \leq \sum_{j=1}^{\lfloor n (T \land N) \rfloor}
         \bigl( |x^k_j - x_j| + |x^k_{j-1} - x_{j-1}| \bigr)
         \left| \be_1^\top \bA^{\nt - j} \be_1 \right|,
 \]
 since
 \ $\left| \be_1^\top \bA^{\nt - j} \be_1 \right| \leq \| \bA^{\nt - j} \|
    \leq C_A$.
We obtain measurability of both \ $\Psi^{N,1}_n$ \ and \ $\Psi^{N,2}_n$, \ hence
 we conclude measurability of \ $\Psi^N_n$.
\ The aim of the following discussion is to show that
 there exists \ $C\subset C_{\Psi,(\Psi_n)_{n\in\NN}}$ \ with \ $C\in\cD_\infty$ \
 and \ $P(\cM \in C) = 1$, \ where \ $C_{\Psi,(\Psi_n)_{n\in\NN}}$ \ is defined in Appendix.
We check that \ $C:=\{ f \in \CC(\RR_+, \RR) : f(0) = 0 \}$ \ satisfies the above mentioned
 conditions.
First note that \ $C=\CC(\RR_+, \RR)\cap\pi_0^{-1}(0)$, \ where \ $\pi_0:\DD(\RR_+,\RR)\to\RR$,
 \ $\pi_0(f):=f(0),$ $f\in\DD(\RR_+,\RR)$.
Using that \ $\CC(\RR_+, \RR)$ \ is a measurable subset of \ $\DD(\RR_+, \RR)$ \
 (see, e.g., Ethier and Kurtz \cite[Problem 3.11.25]{EK})
 and that \ $\pi_0$ \ is measurable (see, e.g., Ethier and Kurtz \cite[Proposition 3.7.1]{EK}),
 we have \ $C\in\cD_\infty$.
\ Fix a function \ $f \in \CC(\RR_+, \RR)$ \ and a sequence \ $(f_n)_{n\in\NN}$
 \ in \ $\DD(\RR_+, \RR)$ \ with \ $f_n \lu f$, \ where \ $\lu$ \ is defined in Appendix.
By the definition of \ $\Psi$, \ we get \ $\Psi(f)\in \CC(\RR_+, \RR)$.
\ Further, we can write
 \begin{align*}
  \Psi_n(f_n)(t)
  &= \sum_{j=1}^\nt
       \left( f_n\left( \frac{j}{n} \right) - f_n\left( \frac{j-1}{n} \right)
              + \frac{\mu_\vare}{n} \right) \be_1^\top \bPi_\bA \be_1 \\
  & \phantom{=\:}
     + \sum_{j=1}^\nt
        \left( f_n\left( \frac{j}{n} \right) - f_n\left( \frac{j-1}{n} \right)
               + \frac{\mu_\vare}{n} \right)
        \be_1^\top (\bA^{\nt - j} - \bPi_\bA) \be_1 ,
   \qquad t\in\RR_+.
 \end{align*}
Using \eqref{HELP_EX_CONVER} and the assumption
 \ $\varrho(A) = \alpha_1+\cdots+\alpha_p=1$,
 we get \ $\be_1^\top \bPi_\bA \be_1 = \frac{1}{\varphi^\prime(1)}$ \ and
 \[
   \sum_{j=1}^\nt
    \left( f_n\left( \frac{j}{n} \right) - f_n\left( \frac{j-1}{n} \right)
           + \frac{\mu_\vare}{n} \right)
   = f_n\left( \frac{\nt}{n} \right) - f_n(0) + \frac{\nt}{n} \mu_\vare .
 \]
Thus we have
 \begin{align*}
  |\Psi_n(f_n)(t)-\Psi(f)(t)|
  &\leq \frac{1}{\varphi^\prime(1)}
        \left| f_n\left( \frac{\nt}{n} \right) - f(t) \right|
        + \frac{\mu_\vare}{n\varphi^\prime(1)} + \frac{|f_n(0)|}{\varphi^\prime(1)}
 \\
  & \phantom{=\:}
   + \sum_{j=1}^\nt
      \left( \left| f_n\left( \frac{j}{n} \right) - f_n\left( \frac{j-1}{n} \right) \right|
             + \frac{\mu_\vare}{n} \right)
      \| \bA^{\nt - j} - \bPi_\bA \| .
 \end{align*}
Here for all \ $T>0$ \ and \ $t \in [0,T]$,
 \begin{align*}
  \left| f_n\left( \frac{\nt}{n} \right) - f(t) \right|
  &\leq \left| f_n\left( \frac{\nt}{n} \right)
               -  f\left( \frac{\nt}{n} \right) \right|
      + \left| f\left( \frac{\nt}{n} \right) - f(t) \right| \\
  &\leq \omega_T(f,n^{-1}) + \sup_{t\in[0,T]} |f_n(t) - f(t)| ,
 \end{align*}
 where \ $\omega_T(f, \cdot)$ \ is the modulus of continuity of \ $f$
 \ on \ $[0,T]$, \ and we have \ $\omega_T(f, n^{-1}) \to 0$ \ since \ $f$ \ is
 continuous (see, e.g., Jacod and Shiryaev \cite[Chapter VI, 1.6]{JSH}).
In a similar way, for all \ $j=1,\ldots,\nt$, \
 \[
   \left| f_n\left( \frac{j}{n} \right)
          - f_n \left( \frac{j-1}{n} \right) \right|
   \leq \omega_T(f, n^{-1}) + 2 \sup_{t\in[0,T]} | f_n(t) - f(t) | .
 \]
By \eqref{rate}, since \ $\varrho(\bA)=1$,
 \[
   \sum_{j=1}^{\nt}
    \left\| \bA^{\nt-j} - \Pi_\bA \right\|
   \leq \sum_{j=1}^{\nt} c_\bA r_\bA^{\nt-j}
   \leq \frac{c_\bA}{1 - r_\bA} .
 \]
Further,
 \[
   |f_n(0)| \leq |f_n(0) - f(0)| + |f(0)|
            \leq \sup_{t\in[0,T]} | f_n(t) - f(t) | + |f(0)| .
 \]
Thus we conclude
 \ $C \subset C_{\Psi,(\Psi_n)_{n\in\NN}}.$
Since \ $\cM_0=0$ \ and, by the definition of a strong solution
 (see, e.g., Jacod and Shiryaev \cite[Definition 2.24, Chapter III]{JSH}), \ $\cM$ \ has
 almost sure continuous sample paths, we have \ $\PP(\cM\in C) = 1$.
\ Consequently, by Lemma \ref{Conv2Funct}, we obtain
 \ $\cX^n=\Psi_n(\cM_n) \distr \Psi(\cM)=\cX$ \ as \ $n\to\infty$.
\proofend

\section{Appendix}
 \label{Appendix}

In the proof of Theorem \ref{main} we will extensively use the following facts
 about the first and second order moments of the sequences
 \ $(X_k)_{k \in \ZZ_+}$ \ and \ $(M_k)_{k \in \ZZ_+}$.

\begin{Lem}\label{Moments}
Let \ $(X_k)_{k\geq -p+1}$ \ be an \INARp\ process defined by \eqref{INARp} such that
 \ $X_0=X_{-1}=\cdots=X_{-p+1}=0$ \ and \ $\EE(\vare_1^2)<\infty$.
\ Then, for all \ $k \in\NN$,
 \begin{gather}
  \EE(X_k)
  = \mu_\vare \sum_{\ell=0}^{k-1} \be_1^\top \bA^\ell \be_1 , \label{mean}\\
  \var(X_k)
  = \sigma_\vare^2 \sum_{\ell=0}^{k-1} (\be_1^\top \bA^{\ell} \be_1)^2
    + \mu_\vare
      \sum_{i=1}^p
       \alpha_i(1-\alpha_i)
         \sum_{j=0}^{k-i-1}\sum_{\ell=0}^{j} (\be_1^\top \bA^{k-j-i-1} \be_1)^2
                                             (\be_1^\top \bA^{\ell} \be_1). \label{var}
 \end{gather}
Moreover,
 \begin{align}
  \EE( M_k \mid \cF_{k-1} )
  & = 0 \qquad \text{for \ $k \in \NN$,}
   \label{m} \\[1mm]
  \EE( M_k M_\ell \mid \cF_{\max\{k,\ell\}-1} )
  & = \begin{cases}
       \alpha_1(1-\alpha_1) X_{k-1} + \cdots + \alpha_p(1-\alpha_p) X_{k-p}
       + \sigma_\vare^2
        & \text{if \ $k = \ell$,} \\
       0 & \text{if \ $k \ne \ell$.}
      \end{cases} \label{Mcond}
 \end{align}
Further,
 \begin{align}
  \EE(M_k) & = 0 \qquad \text{for \ $k \in \NN$,} \label{M} \\[1mm]
  \EE( M_k M_\ell )
   & = \begin{cases}
        \alpha_1(1-\alpha_1) \EE(X_{k-1}) + \cdots
        + \alpha_p(1-\alpha_p) \EE(X_{k-p}) + \sigma_\vare^2
         & \text{if \ $k=\ell$,} \\
        0 & \text{if \ $k \ne \ell$.}
       \end{cases} \label{Cov}
 \end{align}
\end{Lem}

\noindent
\textbf{Proof.} \
We have already proved \eqref{mean}, see \eqref{expect_rec}.
The equality \ $M_k = X_k - \EE(X_k \mid \cF_{k-1})$ \ clearly implies \eqref{m}
 and \eqref{M}.
By \eqref{INARp} and \eqref{Mk},
 \begin{align}\label{HELP_M_DECOMP}
  M_k
  = \sum_{j=1}^{X_{k-1}} \big( \xi_{k,1,j} - \EE(\xi_{k,1,j}) \big) + \cdots
    + \sum_{j=1}^{X_{k-p}} \big( \xi_{k,p,j} - \EE(\xi_{k,p,j}) \big)
    + \big( \vare_k - \EE(\vare_k) \big) .
 \end{align}
For all \ $k\in\NN$, \ the random variables
 \ $\big\{\xi_{k,i,j} - \EE(\xi_{k,i,j}) , \, \vare_k - \EE(\vare_k)
          : j \in \NN, \, i \in \{ 1, \dots, p \} \big\}$
 \ are independent of each other, independent of \ $\cF_{k-1}$, \ and have
 zero mean, thus in the case \ $k = \ell$ \ we conclude \eqref{Mcond} and hence
 \eqref{Cov}.
If \ $k < \ell$, \ then
 \ $\EE( M_k M_\ell \mid \cF_{\ell-1} ) = M_k \EE( M_\ell \mid \cF_{\ell-1} ) = 0$
 \ by \eqref{m}, and thus we obtain \eqref{Mcond} and \eqref{Cov} in the case of
 \ $k \ne \ell$.

By \eqref{REGR} and \eqref{mean}, we conclude
 \[
    X_k - \EE(X_k)
    = \sum_{j=1}^k M_j \be_1^\top \bA^{k-j} \be_1 , \qquad k\in\NN.
 \]
Now, by \eqref{Cov}, \eqref{mean},
 \begin{align*}
  \var(X_k)
 & = \sum_{j=1}^k \sum_{\ell=1}^k
     \EE( M_j M_\ell ) \be_1^\top \bA^{k-j} \be_1 \be_1^\top \bA ^ {k-\ell} \be_1
  = \sum_{j=1}^k
     \EE(M_j^2) (\be_1^\top \bA^{k-j} \be_1)^2 \\
 & = \sum_{j=1}^k
     \left(\sum_{i=1}^p \alpha_i(1-\alpha_i)\EE(X_{j-i}) + \sigma_\vare^2\right)
          (\be_1^\top \bA^{k-j} \be_1)^2 \\
 & = \sigma_\vare^2 \sum_{j=1}^k (\be_1^\top \bA^{k-j} \be_1)^2
      + \sum_{i=1}^p \alpha_i(1-\alpha_i)\sum_{j=1}^k  \EE(X_{j-i}) (\be_1^\top \bA^{k-j} \be_1)^2 ,
 \end{align*}
 and hence, using also that \ $\EE(X_0)=\EE(X_{-1})=\cdots=\EE(X_{-p+1})=0$, \ we get
 \begin{align*}
  \var(X_k)
   = \sigma_\vare^2 \sum_{\ell=0}^{k-1} (\be_1^\top \bA^{\ell} \be_1)^2
      + \sum_{i=1}^p \alpha_i(1-\alpha_i) \mu_\vare \sum_{j=i+1}^k \sum_{\ell=0}^{j-i-1}
         (\be_1^\top \bA^{\ell} \be_1) (\be_1^\top \bA^{k-j} \be_1)^2 \\
   = \sigma_\vare^2 \sum_{\ell=0}^{k-1} (\be_1^\top \bA^{\ell} \be_1)^2
      + \mu_\vare \sum_{i=1}^p \alpha_i(1-\alpha_i) \sum_{j=0}^{k-i-1} \sum_{\ell=0}^j
         (\be_1^\top \bA^{\ell} \be_1) (\be_1^\top \bA^{k-j-i-1} \be_1)^2,
 \end{align*}
 which yields \eqref{var}.
\proofend

\begin{Cor}\label{EEX}
Let \ $(X_k)_{k\geq -p+1}$ \ be a primitive \INARp\ process defined by \eqref{INARp}
 such that  \ $\alpha_1+\cdots+\alpha_p=1$ \ (i.e. unstable), \ $X_0=X_{-1}=\cdots=X_{-p+1}=0$
 \ and \ $\EE(\vare_1^2)<\infty$.
\ Then
 \[
   \EE(X_k) = O(k) , \qquad
   \EE(X_k^2) = O(k^2) , \qquad
   \EE(|M_k|) = O(k^{1/2}) .
 \]
\end{Cor}

\noindent
\textbf{Proof.}
By \eqref{mean},
 \[
   \EE(X_k)
   \leq \mu_\vare \sum_{\ell=0}^{k-1} \| \bA^\ell \|
   \leq C_\bA \mu_\vare k ,
 \]
 where
 \begin{equation}\label{C}
  C_\bA
  := \sup_{\ell \in \ZZ_+} \| \bA^\ell \|
  < \infty.
 \end{equation}
Here \ $C_\bA$ \ is finite since, by \eqref{rate}, \ $C_\bA\leq c_\bA + \Vert\bPi_\bA\Vert.$
\ Hence we obtain \ $\EE(X_k) = O(k)$.
\ We remark that \ $\EE(X_k) = O(k)$ \ is in fact an immediate consequence of part (ii) of
 Proposition \ref{Proposition1}.

We have, by Lyapunov's inequality,
 \begin{align*}
  \EE(|M_k|)
  & \leq \sqrt{\EE(M_k^2)}
    = \left(\sum_{i=1}^p\alpha_i(1-\alpha_i) \EE(X_{k-i}) + \sigma_\vare^2\right)^{1/2} \\
  & \leq  \left(\sum_{i=1}^p
               \alpha_i(1-\alpha_i) \EE(X_{k-i})\right)^{1/2} + (\sigma_\vare^2)^{1/2},
 \end{align*}
 hence we obtain \ $\EE(|M_k|) = O(k^{1/2})$ \ from \ $\EE(X_k) = O(k)$.

Thus we get
 \[
   \EE(X_k^2) = \var(X_k) + (\EE(X_k))^2 = O(k^2) .
 \]
Indeed, by \eqref{var} and \eqref{C},
 \begin{align*}
  \var(X_k)
  & \leq \sigma_\vare^2 \sum_{\ell=0}^{k-1} \| \bA^{\ell} \|^2
         + \mu_\vare
           \sum_{i=1}^p
             \alpha_i(1-\alpha_i)
           \sum_{j=0}^{k-i-1}
           \sum_{\ell=0}^j
            \|\bA^{\ell}\|
              \|\bA^{k-j-i-1}\|^2 \\
  & \leq \sigma_\vare^2 C_\bA^2 k
         + C_\bA^3 \mu_\vare \sigma_\alpha^2 k^2 ,
 \end{align*}
 where \ $\sigma_\alpha^2$ \ is defined in Theorem \ref{main}.
Hence we obtain \ $\EE(X_k^2) = O(k^2)$.
\proofend

Next we recall a result about convergence of step processes towards a diffusion process,
 see Isp\'any and Pap \cite[Corollary 2.2]{IspPap}.
This result is used for the proof of convergence \eqref{Conv_M}.

\begin{Thm}\label{Conv2DiffCor}
Let \ $\gamma : \RR_+\times\RR \to \RR$ \ be a continuous function.
Assume that uniqueness in the sense of probability law holds for the SDE
 \begin{equation}\label{SDE}
  \dd \, \cU_t = \gamma (t,\cU_t) \, \dd \cW_t,
  \qquad t\in\RR_+,
 \end{equation}
 with initial value \ $\cU_0=u_0$ \ for all \ $u_0\in\RR$, \ where
 \ $(\cW_t)_{t \in \RR_+}$ \ is a standard Wiener process.
Let \ $(\cU_t)_{t \in \RR_+}$ \ be a solution of \eqref{SDE} with initial value
 \ $\cU_0 = 0$.

For each \ $n\in\NN$, \ let \ $(U^n_k)_{k\in\NN}$ \ be a sequence of random
 variables adapted to a filtration \ $(\cF^n_k)_{k\in\ZZ_+}$.
\ Let
 \[
   \cU^n_t := \sum_{k=1}^{\nt} U^n_k\,, \qquad t\in\RR_+, \quad n\in\NN.
 \]
Suppose \ $\EE \big( (U^n_k)^2 \big) < \infty$ \ and
 \ $\EE \big( U^n_k \mid \cF^n_{k-1} \big) = 0$ \ for all \ $n,k\in\NN$.
\ Suppose that for each \ $T > 0$,
 \begin{enumerate}
  \item [\textup{(i)}]
        $\sup\limits_{t\in[0,T]}
         \left| \sum\limits_{k=1}^{\nt}
                  \EE( (U^n_k)^2 \mid \cF^n_{k-1} )
                 - \int_0^t
                    \gamma(s,\cU^n_s)^2
                    \dd s \right|
         \stoch 0$,\\
  \item [\textup{(ii)}]
        $\sum\limits_{k=1}^{\lfloor nT \rfloor}
          \EE \big( (U^n_k)^2 \bbone_{\{|U^n_k| > \theta\}}
                    \bmid \cF^n_{k-1} \big)
         \stoch 0$
        \ for all \ $\theta>0$,
 \end{enumerate}
 where \ $\stoch$ \ denotes convergence in probability.
Then \ $\cU^n \distr \cU$ \ as \ $n\to\infty$.
\end{Thm}

In fact, this theorem is a corollary of a more general limit theorem, see Isp\'any and Pap
 \cite[Theorem 2.1]{IspPap}.

Now we recall a version of the continuous mapping theorem.

For a function \ $f \in \DD(\RR_+, \RR)$ \ and for a sequence \ $(f_n)_{n\in\NN}$
 \ in \ $\DD(\RR_+, \RR)$, \ we write \ $f_n \lu f$ \ if \ $(f_n)_{n\in\NN}$
 \ converges to \ $f$ \ locally uniformly, i.e., if
 \ $\sup_{t\in[0,T]} |f_n(t) - f(t)| \to 0$ \ as \ $n \to \infty$ \ for all
 \ $T > 0$.
\ For measurable mappings \ $\Phi : \DD(\RR_+, \RR) \to \DD(\RR_+, \RR)$ \ and
 \ $\Phi_n : \DD(\RR_+, \RR) \to \DD(\RR_+, \RR)$, \ $n \in \NN$, \ we will
 denote by \ $C_{\Phi,(\Phi_n)_{n\in\NN}}$ \ the set of all functions
 \ $f \in \CC(\RR_+, \RR)$ \ such that \ $\Phi(f) \in \CC(\RR_+, \RR)$ \ and
 \ $\Phi_n(f_n) \lu \Phi(f)$ \ whenever \ $f_n \lu f$ \ with
 \ $f_n \in \DD(\RR_+, \RR)$, \ $n \in \NN$.

For deriving convergence \eqref{Conv_X} from convergence \eqref{Conv_M} we will need
 the following version of the continuous mapping theorem.

\begin{Lem}\label{Conv2Funct}
Let \ $(\cU_t)_{t\in\RR_+}$ \ and \ $(\cU^n_t)_{t\in\RR_+}$, \ $n\in\NN$, \ be
 stochastic processes with c\`adl\`ag paths such that \ $\cU^n \distr \cU$ \
 as \ $n\to\infty$.
\ Let \ $\Phi : \DD(\RR_+, \RR) \to \DD(\RR_+, \RR)$ \ and
 \ $\Phi_n : \DD(\RR_+, \RR) \to \DD(\RR_+, \RR)$, \ $n \in \NN$, \ be
 measurable mappings such that
 there exists \ $C\subset C_{\Phi,(\Phi_n)_{n\in\NN}}$ \ with
 \ $C\in\cD_\infty$ \ and \ $\PP(\cU \in C) = 1$.
\ Then \ $\Phi_n(\cU^n) \distr \Phi(\cU)$ \ as \ $n\to\infty$.
\end{Lem}

Lemma \ref{Conv2Funct} can be considered as a consequence of Theorem 3.27 in Kallenberg \cite{Kal},
 and we note that a proof of this lemma can also be found in Isp\'any
 and Pap \cite[Lemma 3.1]{IspPap}.

\section*{Acknowledgements}
Figures 1 and 2 were prepared with SAS Enterprise Guide 4.2.
The authors have been supported by the Hungarian Portuguese Intergovernmental S \& T Cooperation
 Programme for 2008--2009 under Grant No.\ PT--07/2007.
M. Barczy and G. Pap have been partially supported by the Hungarian Scientific Research Fund
 under Grant No.\ OTKA T--079128.
M. Isp\'any has been partially supported by T\'AMOP 4.2.1./B-09/1/KONV-2010-0007/IK/IT
 project, which is implemented through the New Hungary Development Plan co--financed by the
 European Social Fund and the European Regional Development Fund.

\end{document}